\documentclass[10pt]{amsart}
\usepackage{graphicx}
\usepackage{pb-diagram}
\usepackage{colordvi}
\usepackage{graphics,eepic}
\usepackage[all]{xy}
\usepackage{mathrsfs}
\usepackage[centertags]{amsmath}
\usepackage{amsfonts}
\usepackage{amssymb}
\usepackage{amsthm}

\newtheorem{thm}{Theorem}
\newtheorem{cor}{Corollary}[section]
\newtheorem{lem}[cor]{Lemma}
\newtheorem{prop}[cor]{Proposition}
\theoremstyle{definition}
\newtheorem{defn}[cor]{Definition}
\theoremstyle{remark}
\newtheorem{rem}[cor]{Remark}
\newtheorem{notation}[cor]{Notation}
\numberwithin{equation}{section}

\newcommand{\pes}[2]{\langle #1,#2\rangle}
\newcommand{\eps}{\varepsilon}
\newcommand{\N}{\mathbb{N}}
\newcommand{\K}{\mathbb{K}}
\newcommand{\C}{\mathbb{C}}
\newcommand{\R}{\mathbb{R}}

\newcommand{\prk}{\mathbb{P}(\mathbb{K}^n)}

\newcommand{\pp}{{\mathbb{P}\big(\mnk\times\K\big)}\times\mathbb{P}(\mathbb{K}^n)}
\newcommand{\ppuno}{\mathbb{P}\big(\mnk\times\K\big)}
\newcommand{\prc}{\mathbb{P}(\mathbb{C}^n)}
\newcommand{\ppc}{\mathbb{P}\big(\mnc\times\C\big)\times\mathbb{P}(\mathbb{C}^n)}
\newcommand{\ppunoc}{\mathbb{P}\big(\mnc\times\C\big)}

\newcommand{\bi}{\begin{itemize}}
\newcommand{\ei}{\end{itemize}}
\newcommand{\bd}{\begin{description}}
\newcommand{\ed}{\end{description}}
\newcommand{\beq}{\begin{equation}}
\newcommand{\eeq}{\end{equation}}
\newcommand{\beqn}{\begin{eqnarray}}
\newcommand{\eeqn}{\end{eqnarray}}
\newcommand{\beqna}{\begin{eqnarray*}}
\newcommand{\eeqna}{\end{eqnarray*}}
\newcommand{\mnk}{\mathbb{K}^{n\times n}}
\newcommand{\mnc}{\mathbb{C}^{n\times n}}
\newcommand{\pmnk}{\mathbb{P}(\mnk)}
\newcommand{\pmnc}{\mathbb{P}(\mnc)}
\newcommand{\pmnr}{\mathbb{P}(\R^{n\times n})}

\newcommand{\V}{\mathcal{V}}
\newcommand{\W}{\mathcal{W}}

\newcommand{\Z}{\mathbb Z}
\def\P{{\mathbb P}}
\begin{document}

\title[Complexity of Eigenvalue Problem]{Complexity of Path--Following Methods for the
Eigenvalue Problem}%
\author[Diego Armentano]{Diego Armentano}%
\address{Facultad de Ciencias,
Universidad de la Rep\'ublica,
Igu\'a 4225,
11400 Montevideo,
Uruguay}
\email{diego@cmat.edu.uy}



\begin{abstract}
A unitarily invariant projective framework is introduced to analyze the complexity of
path--following methods for the eigenvalue problem.
A condition number, and its relation to the distance to ill--posedness, is given.
A Newton map appropriate for this context is defined,  and a version of Smale's
$\gamma$-Theorem is proven.
The main result of this paper bounds the complexity of path--following methods in terms
of the length of the path in the condition metric.
\end{abstract}
\maketitle

 \par{\centering To the memory of Jean-Pierre Dedieu and Mario
 Wschebor\par}

\medskip

\tableofcontents

\section{Introduction and main results}

\subsection{Introduction and background}

In this paper we study the complexity of path--following methods to solve the eigenvalue
problem: 
$$(\lambda I_n - A)v=0, \quad v \neq 0,$$
where $A \in \K^{n \times n}$, $v \in \K^n$, and $\lambda \in \K$.  
(The set $\K$ denotes $\R$ or $\C$, and $n\geq2$.)
Here, the complexity of an algorithm should be understood as the study of the number of arithmetic operations required to pass from the input to the output.

Many algorithms have been used to solve the eigenvalue problem. A naive approach to solve
this problem would be to compute the characteristic polynomial $\chi_A(z)=\det(z I_n-A)$
of $A$ and
then compute (i.e., find approximations of) its zeros. From these zeros, a
correspondent eigenvector can be computed. 
Unfortunately, in some cases, polynomials $\chi_A(z)$ arising in this way may
be ill--conditioned even when the original matrix $A$ isn't, and therefore,  the
numerical stability may be destroyed under this process.
In practice, algorithms for solving the eigenvalue problem
have avoided this naive approach, and moreover, the tradition in numerical
analysis in the last 20 years seems to go in the opposite direction.
That is, in order to solve a polynomial in one variable the standard procedure
is to search for the eigenvalues of the associated companion matrix; see Trefethen--Bau~\cite{Trefethen-Bau} for example.

Most of the algorithms used in practice for solving the eigenvalue problem may be divided
into two classes:
QR methods (including Hessenberg reduction, single or double shift strategy, 
deflation), and Krylov subspace methods;
see Wilkinson~\cite{Wil}, Golub--Van Loan~\cite{gol}, Stewart~\cite{St2}, 
or Watkins~\cite{Wat} for details. These algorithms are known to be stable but, 
surprisingly, the complexity in not well--understood; indeed, for each of these
methods one of the following fundamental questions is still open: 
\begin{enumerate}
\item For which class of matrices there is guaranteed convergence?
 \item Is there a small bound on the average of number of steps, in a given probabilistic model 
on the set of inputs, to obtain a given accuracy on the output?
\end{enumerate}
The two following examples show that such questions are particularly difficult:
\begin{itemize}
\item Rayleigh quotient iteration fails for a non--empty open set of matrices; see
Batterson-Smillie~\cite{bat1, bat2}.
\item The unshifted QR algorithm is convergent for almost every complex matrix. However,
even for the simple choice of \emph{Gaussian Orthogonal Ensemble} as a
probabilistic model, question (2) remains unanswered; see
Deift~\cite{Deift}. (See also Pfrang et al.~\cite{pfrang-deift-Menon} for some statistics about QR and Toda algorithms.)
\end{itemize}

While in practice many numerical methods are available for computing the eigenvalues and
eigenvectors, until now,  a numerically stable algorithm that provides satisfactory answers to questions (1)
and (2) is not available.  
On the one hand we have algorithms for which we can prove low complexity bounds but  appears to be unstable in practice; see Pan~\cite{Pan96} or Renegar~\cite{Renegar}. On the other hand, we have algorithms which are stable (and even 
efficient) in practice but for which we cannot prove satisfactory complexity bounds. 

The main theoretical open question for the complexity of the eigenvalue problem is to 
provide an algorithm to solve this problem which is numerically stable and
works in average polynomial--time. The present paper may be considered as a step forward to
achieve this goal; the underlying algorithmic approach are the so called homotopy methods. 

In the last three decades path--following methods, or homotopy methods, have
been applied  to solve the
eigenvalue problem. The advantage of using homotopy methods for the eigenvalue problem 
lies in the following facts.
\begin{enumerate}
 \item Path--following methods are numerically stable almost by definition. Being
more precise, in order to follow a path of (problems,solutions) it is enough to compute
a sequence of pairs such that each pair is a mere approximation of the path with some prescribed error; ; see, for instance, the recent analysis
in Briquel et al.~\cite{BCPR:12} for the case of complex polynomial systems. 
 \item  The recent success of homotopy methods for attacking Smale's 17th problem (see
Beltr\'an--Pardo~\cite{Beltran_Pardo-2009-1}, B\"urgisser--Cucker~\cite{BurCuck}, and more recently
Ar\-men\-ta\-no--Shub~\cite{Arm_Shub-2013}) brings some hope to reach the main goal
for the complexity of the eigenvalue problem mentioned above.
\item Given a matrix, the output of the homotopy method is a good approximation
of some eigenvalue  (and its corresponding eigenvector) of the given matrix. In this way
we are
avoiding to deal with the problem of computing an eigenvector from some
approximation
of the eigenvalue of the matrix (which will carry an extra cost and accuracy issues; cf. Remark~\ref{rem:pseudo}).
\end{enumerate}

The homotopy method for the eigenvalue problem was first studied by
Chu~\cite{Chu}, when $A$ is a real symmetric matrix. In
Li--Sauer--York~\cite{Li-Sauer-York} and Li--Sauer~\cite{Li-Sauer}
homotopy methods were given for \textit{deficient} polynomial systems,
and in particular, the general eigenvalue problem is considered; see
Li~\cite{Li} for a general discussion. Since then, a substantial amount of
papers dealing with homotopy methods for solving the eigenvalue
problem have been written; see Lui--Keller--Kwok
\cite{Lui-Keller-Kwok} and references therein. Even though these
methods can achieve spectacular results in practice (even faster than
QR in some cases) the complexity is still an open problem.

\medskip

In this paper we consider the eigenvalue problem as a bilinear
polynomial system of equations and we study the complexity of homotopy
methods to solve it.
Briefly, homotopy methods can be described as follows. The system
$(\lambda I_n - A)v=0,$ $v \ne 0,$ is the endpoint of a path of
problems
$$(\lambda(t) I_n - A(t))v(t)=0, \ v(t) \ne 0, \ 0 \leq t \leq 1,$$
with $(A(1), \lambda(1), v(1)) = (A, \lambda, v)$. Starting from a known triple $(A(0),
\lambda(0), v(0))$ we ``follow'' this path to reach the target system
$(\lambda I_n - A)v=0$. The algorithmic way to do so is to construct a finite number of
triples
$$(A_k, \lambda_k, v_k), \quad 0 \leq k \leq K,$$
with $A_k = A(t_k)$, and $0 = t_0 < t_1 < \cdots < t_K = 1$, and where $(\lambda_k,v_k)$
are approximations of $(\lambda(t_k),v(t_k))$.
The complexity of the algorithm just described (defined more precisely below) is reduced to 
the number $K$ of steps sufficient to validate the approximation, since the arithmetic cost of each iteration is linear in $n$.

The main result of this paper is to relate $K$ with
a geometric invariant, namely, the {\bf condition length} of the path. In the next
paragraphs of this section we give a succinct description of the main definitions in order to state the main results, pointing on our way to the location in this paper when notions are dealt at greater length.

\medskip

\subsection{Solution Variety}\label{ssec:SV}

We begin with the geometric framework of our problem. Since the eigenvalue problem is homogeneous in $v \in \K^n$ and in
$(A, \lambda) \in \K^{n \times n} \times \K$, we define the \textit{solution variety} as
$$ 
\V := \left\{ (A, \lambda, v) \in \P \left( \K^{n \times n} \times \K \right) \times \P \left( \K^{n}\right): \ (\lambda I_n - A)v=0 \right\},
$$
where $\P(\mathbb{E})$ denotes the projective space associated with the vector space $\mathbb{E}$.
We speak interchangeably of a nonzero vector and its corresponding
class in the projective space.

The solution variety $\V$ plays a crucial role in this paper.  
It is a connected smooth manifold of the same dimension as
$\pmnk$.
It is possible to define a natural projection
\begin{equation*}
\pi:\V\to\pmnk\quad\mbox{given by}\quad\pi(A,\lambda,v)=A.
\end{equation*}
This projection, for almost every $(A,\lambda,v)\in\V$, has
a branch of  the inverse image of $\pi$ taking $A\in\pmnk$ to
$(A,\lambda,v)\in\V$. 
This branch of $\pi^{-1}$ is usually called the \emph{input--output map}. In this fashion, we may think of $\pmnk$ as the space of inputs and $\V$ as
the space of outputs. (Section~\ref{subsec:SolVar} provides a detailed exposition of
these facts.)

\medskip

\subsection{Newton's method}\label{sec:BHNMeth}

Given a nonzero matrix $A\in\mnk$, we define the evaluation map $F_A:\K\times\K^n \to \K^n$, by
$$
F_A(\lambda,v):=(\lambda I_n-A)v.
$$

The  \emph{Newton map} associated to $F_A$, is the map $N_A$ on 
$\K\times(\K^n\setminus\{0\})$ given by
$$
N_{A}(\lambda,v):=(\lambda,v)-\big(DF_A(\lambda,v)|_{\K\times v^{\perp}}\big)^{-1}F_A(\lambda,v),
$$
defined for every $(\lambda,v)$ such that
$DF_A(\lambda,v)|_{\K\times v^{\perp}}$ is invertible, where $DF_A(\lambda,v)$ denotes the derivative of $F_A$ at $(\lambda,v)$.
Here $v^\perp$ is the Hermitian complement of $v$ in $\K^n$.

{\sloppypar 
In Section~\ref{sec:NM} we show that the map $N_A$ is well--defined provided that ${\Pi_{v^\perp}(\lambda
I_n-A)\big|_{v^\perp}}$ is invertible (where
$\Pi_{v^\perp}$ denotes the orthogonal projection of $\K^n$ onto $v^\perp$).}
If this is the case, then the map $N_A$ is given by 
$
N_A(\lambda,v)=(\lambda-\dot\lambda,v-\dot v),
$ where
\begin{equation}\label{eq:NAexp}
\dot v  =\left(\Pi_{v^\perp}(\lambda
I_n-A)\big|_{v^\perp}\right)^{-1}\Pi_{v^\perp}(\lambda I_n-A)v, \quad
\dot\lambda =\frac{\pes{(\lambda I_n-A)(v-\dot v)}{v}}{\pes vv}.
\end{equation}
(Here $\pes\cdot\cdot$ denotes the canonical Hermitian product on $\K^n$.)

Let $A\in\mnk$ be a nonzero matrix, and let $(\lambda_0,v_0)\in\K\times\K^n$, $v_0\neq0$.
We say that the triple $(A,\lambda_0,v_0)$ is an 
\emph{approximate solution} of the eigenvalue problem $(A,\lambda,v)\in \V$, when the sequence $(A,N_A^k(\lambda_0,v_0))$, $k=0,1,2,\ldots$ converges immediately quadratically to the eigentriple $(A,\lambda,v)\in\V$, that is, if the given sequence  satisfies
$$
 d_{\mathbb{P}^2}\left((A,N_A^k(\lambda_0,v_0)),(A,\lambda,v)\right)\leq 
 \left(\frac{1}{2}\right)^{2^k-1} d_{\mathbb{P}^2}\left((A,\lambda_0,v_0),(A,\lambda,v)\right),
 $$
 for all positive integers $k$.
Here $d_{\mathbb{P}^2}(\cdot,\cdot)$ is the induced Riemannian distance on $\pp$; see
Section~\ref{ss:MS}.

\begin{rem}\label{rem:AZwelldefpp}
From the expression (\ref{eq:NAexp}), it is easily seen that, if $N_A(\lambda,v)=(\lambda',v')$ then $N_{\alpha A}(\alpha
\lambda,\beta v)=(\alpha\lambda', \beta v')$, for every nonzero scalars $\alpha$ and $\beta$. Hence, the sequence defined above with starting point  $(A,\lambda_0,v_0)$ and that starting at $(\alpha A, \alpha \lambda_0,\beta v_0)$, define the same sequence on $\pp$. Thus the property of being an approximate solution is well--defined on $\pp$.
\end{rem}
\begin{rem}\label{rem:Nwelldefpp}
 From the last remark we conclude that the Newton map $N_A$ induces a map from $\K\times\prk$ into itself (defined almost everywhere).
\end{rem}

\medskip

\subsection{The predictor--corrector algorithm}\label{sub:PCA}

Let
$\Gamma(t)=(A(t),\lambda(t),v(t))$, ${0\leq   t \leq  1}$,   be a
path of eigentriples in
$(\mnk\setminus\{0_n\})\times\K\times\prk$, i.e., $\Gamma\subset\V$.
 To approximate $\Gamma$ by a finite sequence we use the following
predictor--corrector strategy:
given a mesh $0=t_0<t_1<\cdots < t_K=1$ and a pair
$(\lambda_0,v_0)\in\K\times\prk$, we
define
$$
(\lambda_{k+1},v_{k+1}):=N_{A(t_{k+1})}(\lambda_k,v_k),\quad 0\leq k\leq K-1,
$$
(in case it is defined).
We say that the sequence $(A(t_k),\lambda_k,v_k)$, $0\leq k \leq K$,
\emph{approximates} the path $\Gamma(t)$, $0\leq t\leq 1$, when for any $k=0,\ldots,K$, the triple 
$(A(t_k),\lambda_k,v_k)$ is an approximate solution of
the eigentriple $\Gamma(t_k)\in\V$.
In that case we define the \emph{complexity of the sequence} by $K$.

\subsection{Condition of a triple and condition length}

Let $\W\subset \V$ be the set of \emph{well--posed} problems, that is, the set of triples $(A,\lambda,v)\in\V$ such that the input--output map mentioned in Section~\ref{ssec:SV}, taking $A\in\pmnk$ to $(A,\lambda,v)\in\V$,  is locally defined. The set $\W$ is the open set of triples $(A,\lambda,v)\in\V$ such that $\lambda$ is a simple eigenvalue; see Section~\ref{sec:SolVar} for further details.  Let $\Sigma':=\V \setminus \W$ be the \emph{ill--posed variety}.

When $(A,\lambda,v)\in\W$, the
operator $\Pi_{v^\perp}(\lambda I_n-A)|_{v^{\perp}}$ is invertible; see
Section~\ref{sec:SolVar}. The \emph{condition number} of $(A,\lambda,v)\in\W$ is defined by
\begin{equation*}
\mu(A,\lambda,v):=\max\left\{1,\|A\|_F \left\|\left(\Pi_{v^\perp}(\lambda I_n-A)|_{v^{\perp}}\right)^{-1}\right\|\right\},
\end{equation*}
where $\|\cdot\|_F$ and $\|\cdot\|$ are the Frobenius  and operator norms in the space of
matrices; see Section~\ref{ss:CN}.

Let $\Gamma(t)$, $0\leq t\leq 1$, be an absolutely continuous path in $\W$. We define
its \emph{condition--length} as
\begin{equation*}
\ell_\mu(\Gamma):=\int_0^1{\left\|\dot \Gamma(t)\right\|_{\Gamma(t)}\mu\left(\Gamma(t)\right)\,dt},
\end{equation*}
where $\left\|\dot \Gamma(t)\right\|_{\Gamma(t)}$ is the norm of $\dot \Gamma(t)$ in the
unitarily invariant  Riemannian structure on $\V$; see Section~\ref{ss:MS} and
Section~\ref{sec:condlengh}.

\medskip

\subsection{Main results}

Recall that $\V\subset\pp$. Let $\pi_2:\V\to\prk$ be the restriction to $\V$ of the canonical projection $(A,\lambda,v)\mapsto v$. Let $\V_v\subset\V$ be the inverse image of $v$ under $\pi_2$. 
\begin{thm}[Condition Number Theorem]\label{theorem:CNT}
  For $(A,\lambda,v)\in\W$, we get
$$
\mu(A,\lambda,v)\leq \max\left\{1,\frac{1}{\left(1+\frac{|\lambda|^2}{\|A\|_F^2}\right)^{1/2}}
\frac{1}{\sin(d_{\mathbb{P}^2}\left( (A,\lambda,v),\Sigma'\cap\V_v)\right)}\right\}.
$$
\end{thm}
This theorem is a version of the Condition Number Theorem which relates the condition number to the distance to ill--posed problems. Its proof is given in Section~\ref{sec:CNT}.
\medskip

The main theorem concerning the convergence of Newton's iteration is the following.
\begin{thm}[Approximate Solution Theorem]\label{teo:newton}
 There is a universal constant ${c_0>0}$ with the following property. Let $A\in\mnk$ be a
nonzero matrix, and let $(\lambda,v)$, $(\lambda_0,v_0)$ in $\K\times\prk$. If
$(A,\lambda,v)\in\W$ and
 $$d_{\mathbb{P}^2}\left((A,\lambda_0,v_0),(A,\lambda,v)\right)<\frac{c_0}{\mu(A,\lambda,v)},$$
  then, $(A,\lambda_0,v_0)$ is an approximate solution of  $(A,\lambda,v)$. (One may choose $c_0=0.0739$.)
\end{thm}
Theorem~\ref{teo:newton} is a version of the so called
\emph{Smale's $\gamma$-theorem} (see Blum et al. ~\cite{B-C-S-S}), which gives the size of the
basin of attraction of Newton's method.
Different versions of Smale's $\gamma$-theorem for the  symmetric eigenvalue problem and for the  generalized eigenvalue problem are given in Dedieu~\cite{Dedieu} and Dedieu--Shub~\cite{D-S} respectively.

Theorem~\ref{teo:newton} is the main ingredient to prove complexity
results for path--following methods.  The proof of this theorem is included in Section~\ref{sec:NM}.

\medskip

Following these lines our main result is the following.

\begin{thm}[Main Theorem]\label{teo:main}
There is a universal constant $C>0$ such that for any absolutely continuous path
$\Gamma(t)=(A(t),\lambda(t),v(t))$ in $\W$, $0\leq t\leq1$, (with
$\ell_\mu(\Gamma)<\infty$),  there exists a sequence
$(A(t_0),\lambda_0,v_0),\ldots,(A(t_K),\lambda_K,v_K)$ such that, $t_0=0$, $t_K=1$, the
triple $(A(t_k),\lambda_k,v_k)$ is an approximation of $\Gamma(t_k)$, $0\leq k\leq K$,
and
$$
K\leq C\,\ell_\mu(\Gamma)+1.
$$
(One may choose $C=100$.)
\end{thm}
The proof of {Theorem~\ref{teo:main}} is given in {Section~\ref{sec:proofmain}}.

\begin{rem}
The selection of a good starting triple $(A(0),\lambda_0,v_0)$ is an important issue that is beyond  the scope of this paper. 
Nevertheless, we would like to conclude this section suggesting some simple candidates. 
\begin{enumerate}
 \item[(i)] \emph{Rank one matrices}: 
 If we are thinking in a fixed family of starting points, a natural requirement would be that the condition number on this set is small. 
 The family of triples $(vv^*,\|v\|^2,v)\in\V$, for $v\in\K^n\setminus\{0\}$, satisfy this requirement. 
More precisely, it is easily check that this family is a subset of the set of triples where the condition number reaches its minimum value. Note that $0$ is a multiple eigenvalue of $vv^*$.
 (This example is a version of an example conjectured by Shub--Smale~\cite{BezV} to be a good starting point for linear homotopy in the polynomial system case.)
\item[(ii)] \emph{Roots of unity}: Another example is to consider $D=\mbox{Diag}(1,\zeta,\ldots, \zeta^{n-1})$, where $\zeta$ is the $n$th primitive root of the unity. 
Then all the associated eigentriples $(D,1,e_{1}),\,(D,\zeta,e_{2}),\ldots (D,\zeta^{n-1},e_n)$ are well--posed, where   $e_1,\ldots,e_n$ denotes the canonical basis of $\K^n$. Their condition numbers are constant equal to $\sqrt{n}/(2\sin(\pi/n))$.  
(This candidate is a version of the system of polynomials considered in B\"urgisser--Cucker~\cite{BurCuck} as starting point.)
\item[(iii)] \emph{Projection on a subspace}: A different approach is to consider the starting point as a function of the input. 
There are many different strategies to pursue. 
For instance, for a given matrix $A=(a_{ij})_{i,j=1,\ldots,n}\in\mnk$ let
$$
A(0)=\begin{pmatrix}
      a_{11}& a_{12}&\ldots& a_{1n}\\
      0 &a_{22}&\cdots &a_{2n}\\
      \vdots &\vdots& &\vdots\\
      0 &a_{n2} &\cdots& a_{nn}
     \end{pmatrix}.
$$
Then one can consider $(A(0),a_{11},e_1)\in\V$ as starting point. Note that $A(0)$ is the orthogonal projection of the matrix $A$ into the subspace of matrices such that $e_1$ is an eigenvector. This kind of procedures have many advantages compared with the preceding examples. However, the main drawback is that in some cases we cannot assure that $(A(0),a_{11},e_1)$ is well--posed.
(This particular selection of starting point is a version of a starting point considered in Armentano--Shub \cite{Arm_Shub-2013} for the polynomial system case.)   
\end{enumerate}

\end{rem}

\medskip

\subsection{Comments}\label{sec:comments}

In their seminal paper~\cite{BezI}, Shub and Smale relate, in the context of polynomial system solving, the complexity $K$ to three ingredients:
the degree of the considered system, the length of the path $\Gamma(t)$, $0\leq t\leq 1$,  and the condition number of the path.
Precisely, they obtain the complexity
\begin{equation}\label{ineq:complSSBezI}
K\leq C D^{3/2}\ell(\Gamma)\mu(\Gamma)^2,
\end{equation}
where $C$ is a universal constant, $D$ is the maximum of the degrees of the underlying system of polynomials, $\ell(\Gamma)$ is the length of $\Gamma$ in the associated Riemannian structure, and $\mu(\Gamma)=\sup_{a\leq t\leq b} \mu\left(\Gamma(t) \right)$.

In Shub~\cite{BezVI} the complexity $K$ of path--following methods for the polynomial
system solving
problem is analyzed in terms of the condition length of the path.

In the context of polynomial system solving, the eigenvalue problem $Av=\lambda v$ for a
$n\times n$ matrix $A$, with unknowns $\lambda$ and $v$, may be considered as a quadratic
system of equations. When $\K=\C$, by B\'ezout's theorem, after homogenization, one
expects $2^n$ roots. However, this system has at most $n$ isolated roots. Therefore the
eigenvalue problem as a quadratic system belongs to the subset of ill--posed problems, and
hence~\cite{BezI} and~\cite{BezVI} do not apply. For this reason, in order to analyze the
complexity of the eigenvalue problem, a different framework is required.

In Shub--Smale~\cite{BezIV} a unitarily invariant geometric framework is introduced to study the eigenvalue problem, where the input space is the space of matrices $\C^{n\times n}$, and the space of outputs is $\C\times\prc$. In Dedieu--Shub~\cite{D-S} a complexity bound of type (\ref{ineq:complSSBezI}) is obtained for general multi--homogeneous analytic functions, which applies in particular for the generalized eigenvalue problem.

In this paper we pursue a different approach, considering the eigenvalue problem
as a bilinear problem; see Section~\ref{subsection:MD}. 
The main difference of this projective framework compared to the frameworks
mentioned above is that the complexity of the eigenvalue problem is not only unitarily invariant but invariant
under the scaling of the matrix as well, and hence the natural space for the
input is the projective space $\pmnk$.  This approach was greatly inspired by Michael Shub.

\begin{rem}\label{rem:charpolyapproach}
As it was mentioned before, there is a natural connection between the
eigenvalue problem and the problem of finding a root of a polynomial in one
variable. Given a $n\times n$
matrix $A$, the roots of the characteristic polynomial $\chi_A(z)$ are exactly
the eigenvalues of $A$. Therefore in the case $\K=\C$, one may consider this
approach to analyze the complexity of the eigenvalue problem, where proven
 average polynomial--time complexity --for homotopy methods-- are given with
respect to some natural
Gaussian measure on the space of polynomials; see
Armentano--Shub~\cite{Arm_Shub-2013}. However, the push--forward measure
(induced by the map $A\mapsto \chi_A$) of the natural Gaussian measure on
$\C^{n\times n}$ is different from the restriction of the Gaussian measure,
mentioned above, to the space of monic polynomials. Hence, it is not clear how
the complexity of the eigenvalue problem is related to the complexity of
solving polynomials in one variable. In addition, this approach has some
important drawbacks:
\begin{itemize}
 \item The complexity analysis of solving the characteristic
polynomial $\chi_A(z)$ is not invariant under scaling of $A$. More precisely,
the complexity theory for polynomial systems mentioned above, applied to
non--homogeneous polynomials in one variable, is not invariant under scaling of
the roots. This is in contrast to the natural scaling of the eigenvalue and the matrix.
\item The complexity of finding eigenvectors is not considered under this
approach.
\item  The map  $A\mapsto \chi_A$ magnifies the condition number in some cases and
hence the complexity may growth. 
\end{itemize}
On the other hand, the magnifying effect of the condition number mentioned
above is not true in general; see example (e) in Section~\ref{sec:examples}. In particular, if one is interested only in the
complexity of finding eigenvalues (and not eigenvectors),  there exist a possibility that this
approach may improve the complexity in some cases (in contrast to the belief prevailing
among numerical analysts; cf. Trefethen-Bau~\cite{Trefethen-Bau}).
It is hoped that the present paper may help to analyze this issue rigorously and may
provide the elements to give theoretical proofs of the experimental tendencies of
numerical analysts, in particular, the connection of polynomials
in one variable and the companion matrices mentioned at the beginning of this paper. (See Beltr\'an--Shub~\cite{Beltran-Shub} for a similar
discussion and some interesting questions.)
\end{rem}

\begin{rem}
A drawback of homotopy methods is that it is not always possible to lift a path of matrices $A(t)$, $0\leq t\leq 1$,  to the solution variety $\V$. However, this is the only requirement to analyze the complexity. One is able to lift this path when the ``continued'' eigenvalue of the homotopy $A(t)$ remains a simple eigenvalue of $A(t)$ for all $t\in[0,1]$. (A completely different situation can occur when we restrict ourselves to the case $\K=\R$, since the projection $\pi:\V\to\pmnr$ is not even surjective; see Section~\ref{sec:SVconn} for discussion on the connectivity of $\W$ for this case.) The \emph{discriminant variety} $\Sigma:=\pi(\Sigma')$ is an algebraic variety of $\pmnk$ (see Remark~\ref{rem:sigmaprima}), hence when $\K=\C$ it has complex codimension one, and thus almost all straight line paths (at least) can be lifted.

\end{rem}

\begin{rem}\label{rem:pseudo}
Another possible strength of this paper is that one
allows $A$ to be non--normal. Moreover, one of the important arguments for
pseudo--spectral
techniques in numerical linear algebra has been that when applying the QR algorithm one
will only get an approximation to the Schur
form of the matrix, and hence one is solving a perturbed problem. This in turn suggests that
the best error bound one can get is from bounding the perturbation and then one has to resort
to pseudo--spectral theory in order to get a precise and reliable bound. As the
pseudo--spectrum
can be rather wild one may end up with very crude and rather poor error bounds. This is not
the case in our Theorem~\ref{teo:newton}. Thus, if $\ell_\mu(\Gamma)$ can be estimated (and it must not be too large) one may have
a good alternative to the QR method for non--normal problems, and the method would come
with nice error bounds.
\end{rem}

\begin{rem}\label{rem:shortpath}
Armentano~\cite{Arm_tesis-2012} addressed, for the case $\K=\C$,  the
problem of finding short paths for the condition length. It is proved that for
every problem ${(A,\lambda,v)\in\W}$ there exist a path $\Gamma$ in $\W$ joining
$(A,\lambda,v)$ with $(e_1e_1^*,1,e_1)$ (where $e_1$ is the first element of the
canonical basis of $\C^n$) such that
$$
\ell_\mu(\Gamma)\leq C\sqrt{n}\left(C'+\log (\mu(A,\lambda,v))\right),
$$
for some universal constants $C$ and $C'$. This type of results shed some light on the
contribution of our paper. More precisely, this result combined with
Theorem~\ref{teo:main} means that it may be possible to solve the eigenvalue problem with
a small complexity, precisely, logarithmic in the condition number of the ending triple. This motivates the study of short paths or geodesics in the condition metric. Any result on this matter is encouraging and a real challenge.

\end{rem}

\begin{rem}
 Theorem~\ref{teo:main} states the existence of a sequence which approximates $\Gamma\subset\W$ with the given complexity. The sequence is described in the proof of this theorem but is not constructive.  
Our next objective is to transfer these theoretical considerations into a practical algorithm.
This issue will be considered in another paper, and the construction of the path--following algorithm should be analogous to the constructions given by Beltr\'an~\cite{Beltran}, B\"urgisser--Cucker ~\cite{BurCuck}, or Dedieu--Malajovich--Shub~\cite{DedMalShub_Ad}, for the polynomial system case.
\end{rem}

\begin{rem}
For the purpose of this paper, we only require that the derivative of the path
$\Gamma$ is almost everywhere defined, and also that the length and the
condition length of $\Gamma$ are defined. For this reason we ask $\Gamma$ to be
an absolutely continuous path. This is in contrast to Shub~\cite{BezVI} where
the requirement for the path $\Gamma$ is to be a $C^1$ path. The $C^1$
hypothesis seems to be more natural for the implementation of the algorithm.
\end{rem}

\noindent \textit{Note:}  Throughout this paper we work with $\K=\C$. However most
definitions and results can be extended immediately to the case $\K=\R$. Whenever it is
necessary  we shall  state the difference.

\bigskip

\thanks{\textbf{Acknowledgements.}
I am very grateful to  Carlos Beltr\'an, Jean--Pierre Dedieu and Mike Shub for many
useful discussions and continuous support. Specially to Mike Shub to whom I am greatly indebted for proposing
this problem to me and for the uncountable stimulating conversations. 
I also would like to thank Felipe Cucker for the critiques, and
Teresa Krick for help with the understanding of the algebraic properties of the solution
variety.
I  gratefully acknowledge a large number of  comments and remarks by the anonymous
referees which lead me to enhance the presentation of this paper. (In particular Remark~\ref{rem:pseudo} and Remark~\ref{rem:sigmaprima} are due to them.)
Part of this work was elaborated  during the special semester on Foundations of
Computational Mathematics in 2009, in the Fields Institute at Toronto. 
Thanks to the Fields Institute for its hospitality and financial support.\\ 
This paper is part of my PhD thesis, written under the supervision of Jean--Pierre Dedieu (Institut de Math\`ematiques de Toulouse, France) and Mario Wschebor (Universidad de la Rep\'ublica, Uruguay).\\
This research was partially supported by Comisi\'on Sectorial de Investigaci\'on Cien\-t\'i\-fica
(CSIC), Agencia Nacional de Investigaci\'on e Innovaci\'on (ANII), and partially funded by
the Mathamsud grant ``Complexity''.}

\textbf{Dedications.} This paper is in the memory of my beloved friends and
advisors Jean-Pierre Dedieu and Mario Wschebor for their inspiration and dedication. 
Jean-Pierre and Mario passed away, before the defense of my PhD thesis and during the
submission of this paper, after a long battle with cancer. To them all my gratitude.

\section{Solution Variety}\label{sec:SolVar}

\subsection{Introduction}

We start this section defining the canonical metric structures. Following
this we define the solution variety $\V$ and the varieties $\Sigma'$ and
$\Sigma$, and we study some basic geometrical, topological and algebraic properties
of these varieties. 

\subsubsection{Canonical metric structures}\label{ss:MS}

The space  $\K^n$ is equipped with the canonical Hermitian inner product $\pes\cdot\cdot$. The space $\mnk$ is equipped with the Frobenius Hermitian inner product
$$
\pes{A}{B}_F:=\mbox{trace }(B^*A),
$$
where $B^*$ denotes the adjoint of $B$.

In general, if $\mathbb{E}$ is a finite dimensional vector space over $\K$ with the Hermitian inner product $\pes\cdot\cdot$,
we can define an Hermitian structure on $\mathbb{P}(\mathbb{E})$ in the following way:
for $x\in\mathbb{E}$,
$$
\pes{w}{w'}_x:=\frac{\pes{w}{w'}}{\|x\|^2},
$$
for all $w,\,w'$ in the Hermitian complement $x^\perp$ of $x$ in
$\mathbb{E}$, which is a natural representation of the tangent
space $T_x \mathbb{P}(\mathbb{E})$.
Let $d_\mathbb{P}(x,y)$ be the angle between the  vectors $x$ and $y$.

In this way, the space $\pp$ inherits the Hermitian product structure
\beq\label{def:metpp}
\pes{(\dot A,\dot \lambda,\dot v)}{(\dot B,\dot \eta,\dot w)}_{(A,\lambda,v)}:= \pes{(\dot
A,\dot\lambda)}{(\dot B,\dot\eta)}_{(A,\lambda)} + \pes {\dot v}{\dot w}_v,
\eeq
for all $(\dot A,\dot \lambda,\dot v),\;(\dot B,\dot \eta,\dot w)\in (A,\lambda)^\perp\times v^\perp$.

We denote by $d_{\mathbb{P}^2}(\cdot,\cdot)$ the induced Riemannian distance on ${\pp}$.

Throughout this paper we denote by the same symbol $d_\mathbb{P}$ distances on $\prk$, $\pmnk$ and $\mathbb{P}(\mnk\times\K)$.

\medskip

\subsection{The Varieties $\V$, $\Sigma'$ and $\Sigma$}\label{subsec:SolVar}

\begin{defn}
 We define the \textit{solution variety} as
 $$ \V := \left\{ (A, \lambda, v) \in \P \left( \K^{n \times n} \times \K
\right) \times \P \left( \K^{n}\right): \ (\lambda I_n - A)v=0 \right\}.$$
\end{defn}

The solution variety $\V$ is the set of equivalence classes of the set of solutions of $F=0$, where 
$F$ is the bilinear system given by 
\begin{equation}\label{def:F}
F:(\mnk\setminus\{0_n\})\times\K\times(\K^n\setminus\{0\})\to \K^n,\quad F(A,\lambda,v)=(\lambda I_n-A)v. 
\end{equation}
Note that $F(\alpha A,\alpha\lambda,\beta v)=\alpha\beta F(A,\lambda,v)$, for all nonzero
scalars $\alpha$ and $\beta$.
 Therefore  $\V$ is an algebraic subvariety of the product $\pp$. Moreover, since $0$ is a regular value of $F$ we conclude that $\V$ is also a smooth submanifold of $\pp$. 
Its dimension over $\K$ is given by
$$\dim \V=\dim (\mnk\times\K\times\K^n)-n-2=n^2-1.$$
Thus we have concluded the following result.
\begin{prop}
 The solution variety $\V$ is a smooth submanifold of $\pp$ with the same dimension as $\pmnk$. The tangent space $T_{(A,\lambda,v)}\mathcal{V}$ to $\V$ at $(A,\lambda,v)$ is the set of triples $(\dot{A},\dot{\lambda},\dot{v})\in \mnk\times \K\times \K^n$, satisfying
\beq\label{tanspaceV}
 (\dot\lambda I_n-\dot A) v +(\lambda I_n-A)\dot v=0; \quad \pes{\dot A}{A}_F+\dot\lambda\overline{\lambda}=0;\quad \pes{\dot v}{v}=0.
\qed
 \eeq
\end{prop}

\begin{rem}\label{rem:HerStrucV}
{\sloppypar
The solution variety $\V$ inherits the Hermitian structure from $\pp$ defined in (\ref{def:metpp}).}
\end{rem}

We denote by $\pi_1$ and $\pi_2$ the restriction to $\V$ of the canonical projections onto $\mathbb{P}\big(\mnk\times\K\big)$ and $\prk$ respectively.

Note that $\pi_1(\V)\subset\ppuno$ does not include the pair $(0_n,1)$. Therefore we can define the map 
$$
\pi:\V \to \pmnk, \quad \pi:=p\circ \pi_1,
$$
where $p$ is the canonical projection
\begin{equation}\label{def:p}
p:(\ppuno \setminus \{(0_n,1)\})\to\pmnk,\quad p(A,\lambda)=A;
\end{equation}
see the following diagram.
$$
\begin{diagram}
\node{} \node{\underset{(A,\lambda,v)}{\V}} \arrow[2]{s,l}{\pi} \arrow{sw,t}{{{\pi_1}}}   \arrow{se,t}{\pi_2}   \\
\node{\underset{(A,\lambda)}{\mathbb{P}\big(\mnk\times\K\big)\setminus\{(0_n,1)\}}}  \arrow{se,t}{p}  \node{} \node{\underset{v}{\prk}}\\
\node{}\node{\underset{A}{\mathbb{P}(\mnk)}}  \node{}
\end{diagram}
$$

The derivative
\beq\label{def:DPi}
D\pi (A,\lambda,v):T_{(A,\lambda,v)}\mathcal{V} \to T_A \mathbb{P}(\mnk),
\eeq
is a linear operator between spaces of equal dimension.

\medskip

\begin{defn}
We say that the triple $(A,\lambda,v)\in\V$ is  \emph{well--posed} when $D\pi (A,\lambda,v)$ is an isomorphism.
Let $\W$ be the set of well--posed triples, and $\Sigma':=\V \setminus \W$ be the \emph{ill--posed variety}.
Let $\Sigma=\pi(\Sigma')\subset\pmnk$ be the \emph{discriminant variety}, 
i.e., the subset of ill--posed inputs.
\end{defn}

\begin{lem}\label{lem:wp}
The ill--posed variety $\Sigma'$ is the set of triples $(A,\lambda,v)\in \V$ such that $\lambda$ is an eigenvalue of $A$ of algebraic multiplicity $\geq$ 2.
\end{lem}
\begin{proof}
 The linear operator (\ref{def:DPi}) is given by
\begin{equation}\label{eq:projp}
D\pi (A,\lambda,v)(\dot A,\dot\lambda,\dot v)= \dot A+ \frac{\dot \lambda\,\overline\lambda}{\|A\|_F^2} A,
\qquad (\dot A,\dot\lambda,\dot v)\in T_{(A,\lambda,v)}\mathcal{V}.
\end{equation}
According to  (\ref{tanspaceV}), a non--trivial triple in the kernel of $D\pi
(A,\lambda,v)$ has the form
$(\frac{-\dot \lambda\,\overline\lambda}{\|A\|_F^2}\, A,\dot\lambda,\dot v)$, where $\pes{\dot v}{v}=0$, $\dot v\neq 0$, and
$$
\dot\lambda\left( 1+\frac{|\lambda|^2}{\|A\|_F^2}\right) v+ (\lambda I_n-A)\dot v=0.
$$
Then, $\mbox{rank}[(\lambda I_n-A)^2]<n-1$, namely, $\lambda$ is not a simple eigenvalue of $A$. \\
Reciprocally, if the algebraic multiplicity of $\lambda$ is $\geq2$, then there exists
$0\neq w\in v^\perp$ such that $(\lambda I_n-A)w=\alpha v$, for some $\alpha\in\K$. Then,
$(\dot A,\dot\lambda,\dot v)$ given by $\dot
A=\frac{\alpha\overline{\lambda}}{\|A\|_F^2}A$, $\dot\lambda=-\alpha$ and $\dot
v=(1+\frac{|\lambda|^2}{\|A\|_F^2})w$, is a non--trivial triple belonging to $\ker
D\pi(A,\lambda,v)$, and therefore $(A,\lambda,v)\in\Sigma'$.
\end{proof}

\begin{rem}\label{rem:sigmaprima}
From Lemma~\ref{lem:wp} we conclude that $\Sigma'$ is an algebraic subvariety of $\V$.
Since $\Sigma$ is the set of matrices $A$ such that the resultant of $\chi_A(x)$ and
$\chi_A(x)'$ is zero, thus it is an algebraic variety of $\pmnk$; see for example Blum et
al.~\cite{B-C-S-S}.
\end{rem}

\begin{lem}\label{lem:defW} 
One has,
$$
\W=\{(A,\lambda,v)\in\V:\,\Pi_{v^\perp}(\lambda I_n -A)|_{v^\perp}\,\mbox{is invertible}  \}.
$$
\end{lem}
\begin{proof}
Let $(A,\lambda,v)\in\W$. Then, from Lemma~\ref{lem:wp}, $(\lambda I_n-A)v=0$ and
 the rank of ${(\lambda I_n-A)^2}$ is $n-1$. That is, $(\lambda I_n-A)v=0$ and the range of the 
linear operator $(\lambda I_n-A)|_{v^\perp}$, which is a $(n-1)$-dimensional subspace of $\K^n$, does not contain the vector $v$. Thus 
$\Pi_{v^\perp}(\lambda I_n -A)|_{v^\perp}$ is an invertible operator from  $v^\perp$
into itself. Reciprocally, by contradiction, if  we assume that    $(A,\lambda,v)\in\V$
and $\lambda$ is not a simple eigenvalue, then there exists $w\in v^\perp$, $w\neq0$, such
that $(\lambda I_n-A)w=\alpha v$ (for some $\alpha\in\K$).
Then, the linear  operator $\Pi_{v^\perp}(\lambda I_n -A)|_{v^\perp}$ has non--trivial
kernel.
\end{proof}

\subsection{Unitarily invariance}\label{sec:unitaryaction}

Let $\mathbb{U}_n(\K)$ stands for the unitary group when $\K=\C$ or the orthogonal group when $\K=\R$.
The group $\mathbb{U}_n(\K)$ acts on $\prk$ in the natural way. In addition, $\mathbb{U}_n(\K)$ acts on $\pmnk$ by conjugation (i.e., $U\cdot A=UAU^{-1}$), and acts on $\ppuno$ by $U\cdot (A,\lambda)= (UAU^{-1},\lambda)$. These actions define an action on the product space $\pp$, namely, 
\begin{equation}\label{eq:UactionPP}
U\cdot(A,\lambda,v)\mapsto (UAU^{-1},\lambda,Uv),\quad U\in\mathbb{U}_n(\K).
\end{equation}

If $(A,\lambda,v)\in\V$, then $(UAU^{-1},\lambda,Uv)\in\V$. Thus $\V$ is invariant under
the product action (\ref{eq:UactionPP}). Furthermore, if $(A,\lambda,v)\in\W$, the condition
of $\lambda$ being a simple eigenvalue of $A$ is invariant under the action of the group
$\mathbb{U}_n(\K)$ on $(A,\lambda)$, hence  the variety $\W$ is invariant under the action
of $\mathbb{U}_n(\K)$ as well. We have thus proved the following result. 
\begin{lem}\label{lem:unitaryaction}
 The solution variety $\V\subset\pp$, and the subvariety $\W\subset\V$, are invariant under the action of $\mathbb{U}_n(\K)$. 
 \qed
\end{lem}

\begin{rem}\label{rem:unitinv}
The action of the group $\mathbb{U}_n(\K)$ preserves the canonical structures defined on $\prk$, $\pmnk$,  and $\ppuno$. Thus $\mathbb{U}_n(\K)$ acts by isometries on these spaces. In particular, $\mathbb{U}_n(\K)$ acts by isometries on $\V$. In addition, the projections  $\pi_1$, $\pi_2$,  and $\pi$  are $\mathbb{U}_n(\K)$--equivariant, i.e., they commute with the action of $\mathbb{U}_n(\K)$.
\end{rem}

\subsection{Connectivity}\label{sec:SVconn}
In this section we study the connectivity of the varieties $\V$ and $\W$.
\begin{prop} 
The solution variety $\V$ is connected.
\end{prop}
\begin{proof}
Let $\widehat\V\subset(\mnk\setminus\{0_n\})\times\K\times(\K^{n}\setminus\{0\})$ be the inverse image of $\V$ under he canonical quotient projection $(\mnk\setminus\{0_n\})\times\K\times(\K^{n}\setminus\{0\})\to {\pp}$, that is, 
\begin{equation}\label{def:affineV}
{\widehat\V}:=\{(A,\lambda,v)\in(\mnk\setminus\{0_n\}
)\times\K\times(\K^n\setminus\{0\}):\,(\lambda I_n-A)v=0\}.
\end{equation}
It suffices to prove that $\widehat\V$ is connected.

The proof consists in the construction of a path connecting $(A,\lambda,v)\in{\widehat\V}$
 with the triple $(e_1e_1^*-I_n,0,e_1)\in{\widehat\V}$. (Here $e_1$ is the
first element of the canonical basis of $\K^n$, and $e_1^*$ denotes the transpose of the
column vector $e_1$.)

Let $(A,\lambda,v)\in{\widehat\V}$. With out loss of generality we can assume $\|v\|=1$.
Also note that we can connect $v$ with $e_1$ by a rotation path
$\{R_t\}_{t\in[0,1]}\subset\mathbb{U}_n(\K)$. Since $\mathbb{U}_n(\K)$ acts also on
${\widehat\V}$, from Lemma~\ref{lem:unitaryaction}, we can assume $v=e_1$.  \\
Let $A=\begin{pmatrix}
\lambda & a \\
0& \hat{A}
\end{pmatrix}$ be the matrix expression of $A$ in the canonical basis $e_1,\ldots,e_n$ of
$\K^n$, where $a\in\K^{1\times(n-1)}$ and $\hat{A}\in\K^{(n-1)\times(n-1)}$. First assume
that $A\neq \lambda I_n$. Then we can connect $(A,\lambda,e_1)\in{\widehat\V}$ with the
triple $(A-\lambda I_n,0,e_1)\in{\widehat\V}$ by the straight line path $\{(A-t\lambda
I_n,(1-t)\lambda,e_1)\}_{0\leq t\leq 1}\subset{\widehat\V}$. (Note that the condition
$A\neq\lambda I_n$ is required.)  In addition, $(A-\lambda I_n)|_{e_1^\perp}=\begin{pmatrix}
 a \\
\hat{A}-\lambda I_{n-1}
\end{pmatrix}$ is a nonzero $n\times(n-1)$ matrix over $\K$. Since $n\geq2$, the set 
$\K^{n\times(n-1)}\setminus\{0\}$ is connected. Then we can join the matrix
$(A-\lambda I_n)|_{e_1^\perp}$ with the matrix 
$\begin{pmatrix}
0\\-I_{n-1}
\end{pmatrix}$ by a path included in $\K^{n\times(n-1)}\setminus\{0\}$, and thus, $A-\lambda I_n$ and $e_1 e_1^*-I_n$ are connected by a path in $\mnk\setminus\{0_n\}$ such that the first column remains invariant.   
Hence, we can join the triple $(A-\lambda I_n,0,e_1)$ with the end point $(e_1e_1^*-I_n,0,e_1)$ by a path in ${\widehat\V}$.
\\
If $A=\lambda I_n=\begin{pmatrix}
\lambda & 0 \\
0& \lambda I_{n-1}
\end{pmatrix}$, then $\lambda\neq 0$, and therefore $(A,\lambda,e_1)\in{\widehat\V}$ is connected by a straight line path in ${\widehat\V}$ to the triple $(\lambda(I_n-e_1e_1^*),0,e_1)\in{\widehat\V}$. Now we are in the condition of the preceding argument.
\end{proof}

\begin{prop}\label{prop:Wconnectivity}
\begin{enumerate}
\item[(i)] When $\K=\C$, $\W$ is connected.
 \item[(ii)] When $\K=\R$ and $n$ odd, $\W$ has two connected components.
  \item[(iii)] When $\K=\R$ and $n$ even, $\W$ is connected.
\end{enumerate}
\end{prop}
Since $\V$ is connected  and $\Sigma'$ is an algebraic subvariety of $\V$ (see
Remark~\ref{rem:sigmaprima}), the assertion (i) of this proposition follows from fact
that a complex algebraic subvariety of $\V$ can not disconnect it; see for example Blum et
al.~\cite[pp. 196]{B-C-S-S}.

For the proof of assertions (ii) and (iii) we need some definitions and a lemma.

Let $e_1,\ldots,e_n$ be the canonical basis of $\R^n$, and let $\mbox{det}(\cdot)$ be the
determinant function. Let  $v\in\R^n$, $v\neq0$. If $L:v^\perp\to v^\perp$ is a linear
operator, then we define its determinant $\mbox{det}_{v^\perp}(L)$ by
$$
\mbox{det}_{v^\perp}(L)=\det(v,Lv_2,\ldots,Lv_n),
$$
where $v,v_2,\ldots,v_n$ is a positive orthonormal basis of $\R^n$, i.e.,
$\det(v,v_2,\ldots,v_n)=1$.

Let $\widehat \V$ be the set defined in (\ref{def:affineV}) and let
 \begin{equation}\label{def:affineW}
 \widehat \W:=\{(A,\lambda,v)\in{\widehat\V}: \Pi_{v^{\perp}}(\lambda
I_n-A)|_{v^\perp}\,\mbox{invertible}  \}.
 \end{equation}

Let $D:{\widehat\W}\to\R$ be the function given by
$$
D(A,\lambda,v)=\mbox{det}_{v^\perp}( \Pi_{v^{\perp}}(\lambda I_n-A)|_{v^\perp}).
$$
This function is the restriction, to ${\widehat\W}$, of a continuous function and
thus continuous. Let ${\widehat\W}^+$ and ${\widehat\W}^-$ be the inverse image,
under $D$, of the rays $(0,+\infty)$ and $(-\infty,0)$ respectively. 
(It is easily seen that ${\widehat\W}^+$ and ${\widehat\W}^-$ are non--empty.)
Then ${\widehat\W}$
is decomposed in the disjoint union of the open sets ${\widehat\W}^+$ and ${\widehat\W}^-$.

\begin{lem}\label{lem:W+con}
 The sets ${\widehat\W}^+$ and ${\widehat\W}^-$ are connected.
\end{lem}
\begin{proof}
Let $S\mathcal{O}(n)$ be the special orthogonal group, that is, the subgroup of
$\mathbb{U}_n(\R)$ of matrices with determinant equal to one. The proof of this lemma is
divided in several claims.

\underline{Claim I:} The map 
 $D:{\widehat\W}\to\R$ is invariant under the action of  $S\mathcal{O}(n)$ on
$\widehat\W$; hence the action of the special orthogonal group on ${\widehat\W}$ leave
${\widehat\W}^+$ and ${\widehat\W}^-$ invariant:
 \\
Let $(A,\lambda,v)\in{\widehat\W}$ and $U\in S\mathcal{O}(n)$. Let $v,v_2,\ldots,v_n$ be a
positive orthonormal basis of $\R^n$. Since $U$ is orthogonal and has determinant one, 
$Uv,Uv_2,\ldots, Uv_n$ is a positive orthonormal basis as well. Note that    $Uv_i \in  
(Uv)^\perp$ and $$
\Pi_{(Uv)^\perp}(\lambda I_n-UAU^{-1})Uv_i =U(\lambda I_n-A)v_i -\alpha_i Uv,
$$
where $\alpha_i=\pes{U(\lambda I_n-A)v_i}{Uv}$.
Then, 
\begin{align*}
D(UAU^{-1},\lambda,Uv)&= \mbox{det}_{(Uv)^\perp}(\Pi_{(Uv)^\perp}(\lambda I_n-UAU^{-1})|_{(Uv)^\perp})\\
&=\det(Uv,U(\lambda I_n-A)v_2,\ldots,U(\lambda I_n-A)v_n)\\
&=\mbox{det}_{v^\perp}(\Pi_{v^\perp}(\lambda I_n-A)|_{v^\perp}).
\end{align*}  
That is, $D(UAU^{-1},\lambda,Uv)=D(A,\lambda,v)$, for every $U\in S\mathcal{O}(n)$, proving the claim.

Let $\hat\pi_2:{\widehat\V}\to\R^n$ be the canonical projection $(A,\lambda,v)\in{\widehat\V}\mapsto v\in\R^n$. Note that ${\hat\pi_2}^{-1}(e_1)\cap {\widehat\W}$ is the subset of triples $(A,\lambda,v)\in{\widehat\W}$ such that $v=e_1$.

\underline{Claim II:}
Each triple in ${\widehat\W}^+$ is connected, by path in ${\widehat\W}^+$, to a triple in ${{\hat\pi_2}^{-1}(e_1)\cap {\widehat\W}^+}$. Similarly, each triple in ${\widehat\W}^-$ is connected, by path in ${\widehat\W}^-$, to a triple in ${\hat\pi_2}^{-1}(e_1)\cap {\widehat\W}^-$  :
\\  
Let $(A,\lambda,v)\in{\widehat\W}^+$. With out loss of generality we may assume $\|v\|=1$. Let $U_t\in S\mathcal{O}(n)$, $0\leq t\leq 1$, be a path in the special orthogonal group
of $\R^n$ such that $U_0=I_n$ and $U_1$ satisfying $U_1(v)=e_1$. Then by
Lemma~\ref{lem:unitaryaction},  the action $U_t\cdot(A,\lambda,v)$, $0\leq t\leq 1$, is a
path in ${\widehat\W}$ connecting the triple $(A,\lambda,v)\in{\widehat\W}$ to the triple
$(U_1 A U_1^{-1},\lambda,e_1)\in{\hat\pi_2}^{-1}(e_1)\cap {\widehat\W}$. In addition, from
Claim I the path is included in a level set of $D$, therefore, if the path started in a
triple belonging to ${\widehat\W}^+$ (respectively ${\widehat\W}^-$), then the path will
remains on ${\widehat\W}^+$ (respectively on ${\widehat\W}^-$). This finishes the proof of
the Claim II.

 If $(A,\lambda,e_1)\in{\hat\pi_2}^{-1}(e_1)\cap {\widehat\W}$, then we may write
$A=\begin{pmatrix}
\lambda & a\\ 
0&           \hat A \end{pmatrix}$, where ${a=\Pi_{e_1}A|_{e_1^\perp}}$ belongs to $\R^{1\times(n-1)}$ and $\hat A\in\R^{(n-1)\times(n-1)}$ invertible.

Let ${\widehat\W}_0$ be the subset of triples $(B,0,e_1)\in{\widehat\W}^+$ such that $e_1$ is also a left eigenvector of $B$ with eigenvalue $0$, that is,
$$
{\widehat\W}_0=\{(B,0,e_1)\in {\widehat\W}: \Pi_{e_1}B|_{e_1^\perp}=0 \},$$
and let ${\widehat\W}_0^+={\widehat\W}_0\cap{\widehat\W}^+$ and
${\widehat\W}_0^-={\widehat\W}_0\cap{\widehat\W}^-$.

\underline{Claim III:} ${\widehat\W}_0^+$ (respectively
${\widehat\W}_0^-$)  is a deformation retract of ${\hat\pi_2}^{-1}(e_1)\cap
{\widehat\W}^+$ (respectively ${\hat\pi_2}^{-1}(e_1)\cap {\widehat\W}^-$), and therefore
the number of connected components is the same:\\
Let us prove that ${\widehat\W}_0^+$  is a deformation retract of ${\hat\pi_2}^{-1}(e_1)\cap {\widehat\W}^+$. (The proof for the set ${\widehat\W}_0^-$  is analogue.) 
Let $(A,\lambda,e_1)\in{\hat\pi_2}^{-1}(e_1)\cap {\widehat\W}^+$, and write $A=\begin{pmatrix}
\lambda & a\\ 
0&           \hat A \end{pmatrix}$. 
Let $(A_t,\lambda_t,e_1)\in{\widehat\V}$, $0\leq t\leq 1$, 
be the path of triples given by
$$
A_t=\begin{pmatrix}
(1-t)\lambda & (1-t)a\\ 
0&            \hat A-t\lambda I_{n-1}\end{pmatrix},\quad \lambda_t=(1-t)\lambda.
$$
Since $A_te_1=(1-t)\lambda e_1$ and $\Pi_{e_1^\perp}(A_t-\lambda_t I_n)|_{e_1^\perp}=\hat A-\lambda I_{n-1}$, we conclude that $(A_t,\lambda_t,e_1)\in {\hat\pi_2}^{-1}(e_1)\cap {\widehat\W}^+$ for every $t\in[0,1]$. Then, the continuous map 
$$
F:[0,1]\times({\hat\pi_2}^{-1}(e_1)\cap {\widehat\W}^+)\to{\hat\pi_2}^{-1}(e_1)\cap {\widehat\W}^+,\quad F(t,(A,\lambda,v))=(A_t,\lambda_t,e_1),
$$  is a deformation retract of the space ${\hat\pi_2}^{-1}(e_1)\cap {\widehat\W}^+$ onto the subspace ${\widehat\W}_0^+$.
\\
\underline{Claim IV:} ${\widehat\W}_0^+$ and ${\widehat\W}_0^-$ are connected:\\
If $(B,0,e_1)\in{\widehat\W}_0^+$ then we may write $B=\begin{pmatrix}
0 & 0\\ 
0&           \hat B \end{pmatrix}$, where $\hat B=\Pi_{e_1^\perp}B|_{e_1^\perp}$ belongs to the space of $(n-1)\times (n-1)$ invertible matrices with positive determinant, namely $\mathbb{G}l_{n-1}(\R)^+$, which is a connected component of the linear group of $\R^{n-1}$. Then, under the identification 
$$
(B,0,e_1)\in{\widehat\W}_0^+\mapsto \Pi_{e_1^\perp}B|_{e_1^\perp}\in\mathbb{G}l_{n-1}(\R)^+,$$
we have that ${\widehat\W}_0^+$ is homeomorphic to $\mathbb{G}l_{n-1}(\R)^+$ and therefore connected. (The proof is analogue for the set ${\widehat\W}_0^-$.) This proves the claim.

Let us finish the proof  of the lemma showing that ${\widehat\W}^+$ is connected. (The proof of the connectivity of
${\widehat\W}^-$ is analogue.)
From Claim II each triple  in ${\widehat\W}^+$ can be
connected to a triple in ${\hat\pi_2}^{-1}(e_1)\cap {\widehat\W}^+$ by a path in ${\widehat\W}^+$.
From Claim III  the set ${\hat\pi_2}^{-1}(e_1)\cap {\widehat\W}^+$ has the same number of
connected components as ${\widehat\W}_0^+$. Then Claim IV finishes the proof.
\end{proof}

We have concluded from this lemma that $\widehat\W$ has two connected components, namely
$\widehat\W^+$ and $\widehat\W^-$.

\begin{proof}[Proof of Proposition \ref{prop:Wconnectivity}]
The variety $\W$ is the quotient space of ${\widehat\W}$ under the equivalence relation of multiplying by nonzero real numbers the coordinates $(A,\lambda)$ and $v$ respectively. 
Let $q:(\mnk\setminus\{0_n\})\times\K\times(\K^{n}\setminus\{0\})\to \pp$  be the canonical quotient projection. In particular $q(\widehat\W)=\W$. 

Let $(A,\lambda,v)\in{\widehat\W}$.
\\
(ii) Let us assume that $n$ is odd. Since
\begin{align}\label{eq:Dalfa}
D(\alpha A,\alpha \lambda,\beta v)&=\mbox{det}_{(\beta v)^\perp}( \Pi_{(\beta v)^{\perp}}(\alpha\lambda I_n-\alpha A)|_{(\beta v)^\perp})\\
 &=(\alpha)^{n-1}\mbox{det}_{v^\perp}( \Pi_{v^{\perp}}(\lambda I_n-A)|_{v^\perp}), \nonumber
\end{align}
for all $\alpha,\,\beta\in\R\setminus\{0\}$, and $n-1$ even, we conclude that the set of
equivalence classes, of the quotient projection $q$, are include in one and only one of the
connected components of $\widehat\W$. Hence $q(\widehat\W^+)\cap
q(\widehat\W^-)=\emptyset$. Furthermore, since the quotient projection $q$ is open and
continuous we conclude, from Lemma \ref{lem:W+con}, that $q(\widehat\W^+)$ and
$q(\widehat\W^-)$ are (non--empty) open and connected sets. Thus $\W$ has two connected components,
namely, $q(\widehat\W^+)$ and $q(\widehat\W^-)$.
\\
(iii) When $n$ is even, then $n-1$ is odd and therefore, from (\ref{eq:Dalfa}), the triples $(A,\lambda,v)$ and
$(-A,-\lambda,v)$
are in different components of ${\widehat\W}$, though they are equivalent triples in $\W$.
Therefore every triple $(A,\lambda,v)\in\W$ has a  representative in ${\widehat\W}^+$.
Hence, from Lemma \ref{lem:W+con}, we obtain that $\W$ is connected.
\end{proof}

\subsection{Multidegree of $\V$}\label{subsection:MD}

The eigenvalue problem as a quadratic system belongs to the subset of ill--posed problems; 
see Section~\ref{sec:comments}. The aim of this section is to prove that the bilinear
approach considered in this paper gives the correct number of roots.

There are many different strategies to prove this with origins in algebraic geometry or algebraic topology. Yet we have not found a convenient proof in the literature to cite. One referee suggested that a possible scheme is to fix the matrix $A$ and homogenize the eigenvalue $\lambda$ in order to obtain a bi--homogeneous system in the variables $(\lambda,v)$, and then apply some basic toric variety theory to obtain the result. Here we pursue a different strategy which make use in a more natural way our bilinear approach in $(A,\lambda)$ and $v$. 

For the sake of simplicity in the exposition we restrict ourself to the case $\K=\C$.
This section follows closely D'Andrea--Krick--Sombra~\cite{DKS}.

Since $\V$ is an algebraic subvariety of the product space $\ppc$, there is a natural
algebraic invariant associated to $\V$, namely, the \emph{multidegree} of $\V$. This
invariant is given by the numbers $\deg_{(n^2-1-i,i)}(\V),\,i=0,\ldots, n-1,$ where
$\deg_{(n^2-1-i,i)}(\V)$ is the number of points of intersection of $\V$ with the product
$\Lambda\times \Lambda'\subset\ppc$, where $\Lambda\subset\ppunoc$ and
$\Lambda'\subset\prc$ are generic $(n^2-1-i)$-codimension plane and $i$-codimension plane
respectively; see Fulton~\cite{Fulton}.

\begin{lem}\label{lem:multidegree}
One has,
 $$
 \deg_{(n^2-1-i,i)}(\V)=\binom{n}{i+1}, \qquad\mbox{for}\quad i=0,\ldots,n-1.
 $$ 
\end{lem}
In order to give a proof of this lemma we recall some definitions from \emph{intersection
theory}; see Fulton~\cite{Fulton}. (See also D'Andrea--Krick--Sombra~\cite{DKS}.)

The Chow ring of $\ppc$ is the graded ring 
$$
\mathcal{A}^*\left(\ppc\right)=\Z[\omega_1,\omega_2]/(\omega_1^{n^2+1},\omega_2^n),
$$
where $\omega_1$ and $\omega_2$ denotes the rational equivalence classes of the inverse images of hyperplanes of $\ppunoc$ and $\prc$, under the projections $\ppc\to \ppunoc$ and $\ppc\to\prc$ respectively.

Given a codimension $n$ algebraic subvariety $\mathcal X\subset\ppc$, the class of $\mathcal X$ in the Chow ring is
$$
[\mathcal X]=\sum_{i=0}^{n-1}\deg_{(n^2-1-i,i)}(\mathcal X)\,\omega_1^{i+1}\omega_2^{n-1-i}\in \mathcal{A}^*\left(\ppc\right). 
$$
\begin{proof}[Proof of Lemma~\ref{lem:multidegree}]
Let $F_i$, $(i=1,\ldots,n)$, be the coordinate functions of $F$ defined in (\ref{def:F}).
Since $F_i$ is bilinear for each $i$,  the class of $\{F_i=0\}$, as a subset of $\ppc$, is given by
$$
[\{F_i=0\}]=\omega_1+\omega_2\in \mathcal{A}^*\left(\ppc\right),\quad (i=1,\ldots,n).
$$
Then, the class of $\V$ in the Chow ring is
$$
[\V]=[\{F_1=0\}\cap\cdots\cap\{F_n=0\}]=\prod_{i=1}^n[\{F_i=0\}],
$$
where the last equality follows from B\'ezout identity.
Therefore one gets
$$
[\V]=(\omega_1+\omega_2)^n \equiv \sum_{\ell=1}^{n}\binom{n}{\ell}\omega_1^{\ell}\omega_2^{n-\ell},
$$
that is, $\deg_{(n^2-1-i,i)}(\V)=\binom{n}{i+1}.$
\end{proof}

From Lemma~\ref{lem:multidegree} we obtain that the number of points of intersection of $\V$ with the product $\Lambda\times\prc$ (for $\Lambda\subset\ppunoc$ a generic hyperplane of codimension $n^2-1$) is $n$. In particular the inverse image of $A\in\pmnc\setminus\Sigma$, under the projection  $\pi:\V\to\pmnc$, is the intersection of $\V$ with $\Lambda_A\times\prc$, where $\Lambda_A$ is a particular hyperplane of codimension  $n^2-1$ in $\ppunoc$, namely, $\Lambda_A$ is the projective line containing the pair of points $\{(A,0),(A,1)\}\in\ppunoc$. However, the family of all projective lines $\Lambda_A$, varying $A$ on $\pmnc$, is not a generic family on $\ppunoc$. In the next proposition we prove that, actually, the family $\Lambda_A$ for $A\in\pmnc\setminus\Sigma$ is included in the generic family of hyperplanes of codimension $n^2-1$ of $\ppunoc$ satisfying Lemma~\ref{lem:multidegree} for $i=0$.
\footnote{If $A\in\pmnc\setminus\Sigma$ then $A$ has $n$ distinct eigenvalues, thus the cardinal number of $\pi^{-1}(A)$ coincides with $\deg_{(n^2-1,0)}(\V)$. In Proposition ~\ref{prop:corrnumroots} we give an independent proof of this fact which we consider interesting \textit{per se}.}

Let $A$ be a finite set. We denote by $\#A$ the cardinal number of $A$.
\begin{prop}\label{prop:corrnumroots}
For all $A\in \pmnc\setminus\Sigma$ we have $\#\pi^{-1}(A)=\deg_{(n^2-1,0)}(\V)$. 
\end{prop}
\begin{proof}
Recall that if $A\in\pmnc\setminus\Sigma$ then the projection $\pi$ is a local
diffeomorphism between a neighbourhood of each inverse image of $A$  in $\V$ and
a neighbourhood of $A$ in $\pmnc\setminus\Sigma$. Then the number of inverse
image is locally constant on $\pmnc\setminus\Sigma$. Furthermore, since $\Sigma$ is an
algebraic subvariety of $\pmnc$ (see Remark~\ref{rem:sigmaprima}), then
$\pmnc\setminus\Sigma$ is connected; cf. proof (i) of
Proposition~\ref{prop:Wconnectivity}. Thus the number of inverse images under $\pi$ is
constant on $\pmnc\setminus\Sigma$.
If $(A,\lambda,v)\in\V\setminus\Sigma'$ then, from Lemma~\ref{lem:wp},  the number of inverse image of $(A,\lambda)$, under ${\pi_1:\V\to\ppunoc}$, is one. Hence the restriction $\pi_1|_{(\V\setminus\Sigma')}:(\V\setminus\Sigma')\to \ppunoc$ is a bijective map onto its image $\pi_1(\V\setminus\Sigma')$. 
Therefore  given $A\in\pmnc\setminus\Sigma$,
we have ${\#\pi^{-1}(A)=\#(p|_{\pi_1(\V)})^{-1}(A)}$, where $p$ is the projection map given in (\ref{def:p}). 
In addition, from~\cite[Corollary 5.6]{Mumford}, we get that $\#(p|_{\pi_1(\V)})^{-1}(A)=\deg \pi_1(\V)$, where $\deg$ is the degree of the projective algebraic subvariety $\pi_1(\V)\subset \ppunoc$. 
Since $\dim\pi_1(\V)=\dim(\V)$ and the fact that $\pi_1|_{(\V\setminus\Sigma')}:(\V\setminus\Sigma')\to\pi_1(\V\setminus\Sigma')$ is bijective, we get that $\#(\Lambda\times\prc)\cap\V=\#\Lambda\cap\pi_1(\V)$, for a generic $(n^2-1)$-codimension plane $\Lambda\subset\ppunoc$.  
Then we conclude that $\deg\pi_1(\V)=\deg_{(n^2-1,0)}(\V)$.
\end{proof}

\begin{rem}
 From Proposition~\ref{prop:corrnumroots} and Lemma~\ref{lem:multidegree} we get that the restriction of the projection ${\pi|_{(\V\setminus\pi^{-1}(\Sigma))}:(\V\setminus\pi^{-1}(\Sigma))\to \pmnc\setminus\Sigma}$ is an $n$-fold covering map. 
\end{rem}

\section{Condition number}\label{ss:CN}

\subsection{Introduction}

In this section we introduce the eigenvalue and eigenvector condition numbers. 
We study some basic properties of these condition numbers and we show some examples. 
We define the condition number of the eigenvalue problem. 
We discuss the condition number theorem for this framework, which relates the condition
number  with the distance to ill--posed problems.
In the last part of this section we study the rate of change of condition numbers.

\subsection{Eigenvalue and eigenvector condition numbers}\label{sec:EigEigCN}

When $(A,\lambda,v)$ belongs to $\W$, according to the implicit function theorem, $\pi$ has an inverse defined in some neighbourhood $\mathcal{U}_A\subset \mathbb{P}(\mnk)$  of $A$ such that $\pi^{-1}(A)=(A,\lambda,v)$.
This map $\displaystyle{\mathscr{S}=\pi^{-1}|_{\mathcal{U}_A}:\mathcal{U}_A\to\V}$ is called the \emph{solution map}. It associates to any matrix ${B\in \mathcal{U}_A}$ the eigentriple $(B,\lambda_B,v_B)$ close to $(A,\lambda,v)$.
Its derivative
\begin{equation*}
D\mathscr{S}(A,\lambda,v):T_A \mathbb{P}(\mnk)\to T_{(A,\lambda,v)}\V,
\end{equation*}
is called the \emph{condition operator} at $(A,\lambda,v)$.

If $(A,\lambda,v)\in \W$, the derivative $D\mathscr{S}(A,\lambda,v)$ associates to each $\dot B\in T_A \mathbb{P}(\mnk)$ a triple 
$(\dot A,\dot\lambda,\dot v)\in T_{(A,\lambda,v)}\V$. This
association  defines two linear maps,
$$
D\mathscr{S}_\lambda(A,\lambda,v)\dot B=(\dot A,\dot\lambda)
\quad\mbox{and}\quad
D\mathscr{S}_v(A,\lambda,v)\dot B=\dot v,
$$
namely, the condition operators of the eigenvalue and eigenvector  respectively.

Recall that $\mathbb{P}(\mnk)$ is equipped with  the canonical Hermitian
structure induced by the Frobenius Hermitian product on $\mnk$.

\begin{defn}\label{def:eigvaleigveccn}
The condition numbers of the eigenvalue and eigenvector, at $(A,\lambda,v)\in\W$, are defined by
\begin{align*}
\mu_\lambda(A,\lambda,v)&:=\sup_{\substack{\dot B\in A^\perp\\ \|\dot
B\|_F=\|A\|_F}} \|D\mathscr{S}_\lambda(A,\lambda,v)\dot B\|_{(A,\lambda)},
\\
\mu_v(A,\lambda,v)&:=\sup_{\substack{\dot B\in A^\perp\\ \|\dot B\|_F=\|A\|_F}}
\|D\mathscr{S}_v(A,\lambda,v)\dot B\|_{v}.
\end{align*}
\end{defn}

\begin{prop}\label{prop:muvmulambdaexp}
Let $(A,\lambda,v)\in\W$. Then,
\begin{enumerate}
\item[(i)]
$$\mu_\lambda(A,\lambda,v)=\frac{1}{1+\frac{|\lambda|^2}{\|A\|_F^2}}\,\left[
1+\frac{\|v\|^2\,\|u\|^2}{|\pes vu|^2} \right]^{1/2},$$
where $u\in\K^n$ is any left eigenvector of $A$ with eigenvalue $\lambda$: a nonzero vector satisfying $(\lambda I_n-A)^*u=0$.

\item[(ii)]
$$\mu_v(A,\lambda,v)=\|A\|_F\,\|(\Pi_{v^\perp}(\lambda
I_n-A)|_{v^\perp})^{-1} \|,$$
where $\|\cdot\|$ is the operator norm.
\end{enumerate}
\end{prop}

\begin{rem}
 ${\Pi_{v^\perp}(\lambda I_n-A)|_{v^\perp}}$ is a linear map from the Hermitian
complement of $v$ in $\K^n$ into itself. Hence the operator norm of its inverse
is independent of the representative of $v$ in $\prk$.
\end{rem}

For the proof of Proposition~\ref{prop:muvmulambdaexp} we need two lemmas.

\begin{lem}\label{lem:condop} Let $(A,\lambda,v)\in\W$. Then for $\dot B\in T_A \mathbb{P}(\mnk)$, one gets:
\begin{enumerate}
 \item[(i)] 
$$
D\mathscr{S}_\lambda(A,\lambda,v)\dot B= \big(\dot B -\dot\lambda \frac{\overline\lambda}{\|A\|_F^2}A,\dot \lambda \big),\,\mbox{ where }\, \dot\lambda=\frac{\pes{\dot B v}{u}}{\left(1+\frac{|\lambda|^2}{\|A\|_F^2}\right)\pes vu},
$$
where $u\in\K^n$ is any left eigenvector of $A$ with eigenvalue $\lambda$;
\item[(ii)]
$$
D\mathscr{S}_v(A,\lambda,v)\dot B= (\Pi_{v^\perp}(\lambda I_n-A)|_{v^\perp})^{-1}\Pi_{v^\perp}(\dot B v).
$$

\end{enumerate}
\end{lem}
\begin{proof}
(i): Let $\dot B\in A^\perp$, and let $(\dot A,\dot \lambda)\in
(A,\lambda)^\perp$ such that $D\mathscr{S}_\lambda(A,\lambda,v)\dot B=(\dot
A,\dot \lambda)$.
Then, by the definition of $D\mathscr{S}_\lambda(A,\lambda,v)$ and
(\ref{eq:projp}) we get 
\beq\label{eq:Bdot}
\dot B = \dot A + \frac{\dot \lambda\,\overline\lambda}{\|A\|_F^2}\,A.
\eeq
Let $u\in\K^n$ be any left eigenvector of $A$ with eigenvalue $\lambda$.
Since $u$ is in the Hermitian complement of the range of  $(\lambda
I_n-A)|_{v^\perp}$, then, from (\ref{tanspaceV}) we get
$\pes{\dot A v}{u}=\dot \lambda \pes{v}{u}$. Furthermore, since
$(A,\lambda,v)\in\W$, then $v\notin \mbox{Im}(A-\lambda I_n)$, and $\pes{v}{u}\neq 0$. Thus
\beq\label{eq:lambdadot}
\dot\lambda =\frac{\pes{\dot A v}{u}}{\pes{v}{u}}.
\eeq
From (\ref{eq:Bdot}) and (\ref{eq:lambdadot}) follows
$$
\dot \lambda= \frac{1}{1+\frac{|\lambda|^2}{\|A\|_F^2}}\,\frac{\pes{\dot B v}{u}}{\pes vu}.
$$
(ii): 
From (\ref{tanspaceV}) again one gets $\Pi_{v^\perp}(\lambda I_n-A)\dot v=\Pi_{v^\perp}\dot A v$. 
Furthermore, since $(A,\lambda,v)\in\W$, then
$$
\dot v= (\Pi_{v^\perp}(\lambda I_n-A)|_{v^\perp})^{-1}\Pi_{v^\perp}( \dot A v).
$$
Since, from (\ref{eq:Bdot}), one has $\Pi_{v^\perp}( \dot B v) =\Pi_{v^\perp}( \dot A v)$,
the result follows.
\end{proof}

\begin{lem}\label{lem:mulambdatec}
 Let $(A,\lambda,v)\in\W$ and let $u\in\K^n$ be any left eigenvector of $A$ with eigenvalue $\lambda$. Then 
$$
 \sup_{\substack{\dot B\in A^\perp\\ \|\dot B\|_F=\|A\|_F}} \left|\pes{\dot B v}{u}\right|=
\|A\|_F\,\sqrt{ \|v\|^2\,\|u\|^2-\frac{|\lambda|^2}{\|A\|_F^2}\,|\pes vu|^2  }.
$$
\end{lem}
\begin{proof}
 Note that $\pes{M v}{u}=\pes{M}{uv^*}_F$ for every matrix $M\in\mnk$.
Write 
\begin{equation*}
 uv^*=\pes{uv^*}{\frac{A}{\|A\|_F}}_F\frac{A}{\|A\|_F}+\alpha C,
 \end{equation*}
  where $C\in A^\perp$ and $\|C\|_F=1$. Then
 $$
  \sup_{\substack{\dot B\in A^\perp\\ \|\dot B\|_F=\|A\|_F}} \left|\pes{\dot B v}{u}\right|=
  \sup_{\substack{\dot B\in A^\perp\\ \|\dot B\|_F=\|A\|_F}} \left|\pes{\dot B}{\alpha C}_F\right|=|\alpha|\,\|A\|_F.
 $$
 Furthermore, 
$$
 |\alpha|^2=\|uv^*\|_F^2-|\pes{uv^*}{\frac{A}{\|A\|_F}}_F |^2,
$$
where $\|uv^*\|_F=\|u\|\,\|v\|$. Since $Av=\lambda v$, then
$|\pes{uv^*}{\frac{A}{\|A\|_F}}_F|=\frac{|\lambda|}{\|A\|_F}|\pes{u}{v}|$. 
\end{proof}

\begin{proof}[Proof of Proposition~\ref{prop:muvmulambdaexp}]
 (i):
From {Lemma~\ref{lem:condop}}, for any $\dot B$ such that $\pes{\dot B}{A}_F=0$, and $\|\dot B\|_A=1$, 
\begin{align}
\|D\mathscr{S}_\lambda(A,\lambda,v)\dot B\|_{(A,\lambda)}^2& = \frac{\|A\|_F^2 +|\dot\lambda|^2
\left(1+\frac{|\lambda|^2}{\|A\|_F^2}\right)}{\|A\|_F^2+|\lambda|^2} 
\label{eq:normadslambda}\\
&=  \frac{\|A\|_F^2 +\left|\frac{\pes{\dot B v}{u}}{\pes
vu}\right|^2
\left(1+\frac{|\lambda|^2}{\|A\|_F^2}\right)^{-1}}{\|A\|_F^2+|\lambda|^2}.
\nonumber
\end{align} 
Then, the proof of (i) can be deduced from Lemma~\ref{lem:mulambdatec}.

(ii): Since $Av=\lambda v$, we have  $\Pi_{v^\perp}( \dot B v) =\Pi_{v^\perp}( (\dot B +\alpha A) v) $, for any $\alpha\in\K$ and $\dot B\in A^{\perp}$. Then, from {Lemma~\ref{lem:condop}} we get:
\begin{align*}
\mu_v(A,\lambda,v)&=\sup_{\substack{\dot B\in A^\perp\\ \|\dot B\|_F=\|A\|_F}} \left\| (\Pi_{v^\perp}(\lambda I_n-A)|_{v^\perp})^{-1}\Pi_{v^\perp}( \dot B v)   \right\|_{v} \\
&= \sup_{\substack{\dot B\in \mnk\\ \|\dot B\|_F=1}} \|A\|_F\, \left\|
(\Pi_{v^\perp}(\lambda I_n-A)|_{v^\perp})^{-1}\Pi_{v^\perp}( \dot B v)
 \right\|_{v} .
\end{align*}
Since $\{\Pi_{v^\perp}( \dot B v):\; \dot B\in \mnk,\, \|\dot B\|_F=1\} $ fills the ball of radius $\|v\|$ in $v^\perp$, the result follows.
\end{proof}

\subsection{Some basic properties}

In the next paragraphs we show some basic properties concerning the condition numbers $\mu_\lambda$ and $\mu_v$.
\begin{prop}\label{prop:mulammuvunitinv}
The condition numbers $\mu_\lambda$ and $\mu_v$ are invariant under the action of the group $\mathbb{U}_n(\K)$, that is, 
\begin{align*}
\mu_\lambda(UAU^{-1},\lambda,Uv)&=\mu_\lambda(A,\lambda,v),\\
 \mu_v(UAU^{-1},\lambda,Uv)&=\mu_v(A,\lambda,v),
\end{align*} 
for every $U\in\mathbb{U}_n(\K)$.
\end{prop}

\begin{rem}
 The proof of Proposition~\ref{prop:mulammuvunitinv} can be deduced from the expressions of $\mu_\lambda$ and $\mu_v$ given in Proposition~\ref{prop:muvmulambdaexp}. However, we prefer to give a different proof which emphasize the fact that the property of the condition numbers of being unitarily invariant resides on the natural election of our Hermitian structures given in our geometric framework. 
\end{rem}
\begin{proof}[Proof of Proposition~\ref{prop:mulammuvunitinv}]
The condition operators of the eigenvalue and eigenvector are given by the derivative of the (locally defined) maps  ${\mathscr{S}_\lambda=\pi_1\circ\mathscr{S}}$ and ${\mathscr{S}_v=\pi_2\circ\mathscr{S}}$ respectively.  From  Remark~\ref{rem:unitinv} the projections $\pi$, $\pi_1$ and $\pi_2$ are $\mathbb{U}_n(\K)$-equivariants, hence the action of $\mathbb{U}_n(\K)$ commutes  with $\mathscr{S}$, $\mathscr{S}_\lambda$ and $\mathscr{S}_v$. In addition, since $\mathbb{U}_n(\K)$ acts by isometries on all the intervening spaces, the result follows.
\end{proof}

\begin{lem}\label{lem:CNlowerbound}
 The condition numbers $\mu_\lambda$ and $\mu_v$ are bounded below by $1/\sqrt{2}$.
\end{lem}
\begin{proof}
 Let $(A,\lambda,v)\in\W$. Since $|\lambda|\leq\|A\|_F$, the proof for $\mu_\lambda$ follows immediately by Proposition~\ref{prop:muvmulambdaexp}. 
 
For the proof for $\mu_v$, first fix a representative of $(A,\lambda,v)\in\W$ such that $\|A\|_F=1$ and $\|v\|=1$. In addition,
since the action of $\mathbb{U}_n(\K)$ on $\prk$ is transitive, by Proposition~\ref{prop:mulammuvunitinv}  we may assume that $v$ is the first element of the
canonical basis. Then $A$ has the form
$\begin{pmatrix}
 \lambda & a\\
0 & \hat A
\end{pmatrix}$, where $a\in\K^{1\times (n-1)}$ and $\hat
A\in\K^{(n-1)\times(n-1)}$. 
Under these assumptions $\Pi_{v^\perp}(\lambda
I_n-A)|_{v^\perp}=\lambda I_{n-1}-\hat A$, hence
$$
\|\Pi_{v^\perp}(\lambda I_n-A)|_{v^\perp} \|\leq
\|\hat A \| +|\lambda|\leq \|\hat A\|_F +|\lambda|\leq \sqrt{2}\|A\|_F,
$$
where last inequality follows from the inequality $x+y\leq\sqrt{2}(x^2+y^2)^{1/2}$, for $x$, $y$, in $\R$.
Since we assume $\|A\|_F=1$, we obtain, $\|\Pi_{v^\perp}(\lambda I_n-A)|_{v^\perp} \|\leq \sqrt{2}$.
Therefore from Proposition~\ref{prop:muvmulambdaexp} we get
$$
1= \|\left( \Pi_{v^\perp}(\lambda I_n-A)|_{v^\perp} \right)^{-1}\,  \Pi_{v^\perp}(\lambda I_n-A)|_{v^\perp} \|\leq\sqrt 2 \mu_v(A,\lambda,v).
$$
\end{proof}
\begin{rem}
 Examples (a) and (b) of Section~\ref{sec:examples} show that the lower bound in Lemma~\ref{lem:CNlowerbound} is sharp. 
\end{rem}

\begin{rem}\label{rem:muvtrans}
Let $(A,\lambda,v)\in\W$. Then $(A+\alpha
I_n,\lambda+\alpha,v)\in\W$, for all $\alpha\in\K$, and $\mu_\lambda (A+\alpha
I_n,\lambda+\alpha,v)$ is constant independent of $\alpha$.
 On the other hand, this is not the case for the eigenvector
condition number. More precisely, 
$$
{\mu_v(A+\alpha
I_n,\lambda+\alpha,v)}= \frac{\|A+\alpha I_n\|_F}{\|A\|_F}{\mu_v(A,\lambda,v)}.
$$
In particular, it is an easy exercise to check that $\mu_v(A+\alpha
I_n,\lambda+\alpha,v)$ is minimized, as a function of $\alpha$, when the matrix
$A+\alpha I_n$ has trace equal to zero, namely, ${\alpha=-\mbox{tr}(A)/n}$. 
As we see in the next section, this procedure may improve drastically $\mu_v$ in some cases, and thus, it could be used as a natural pre--conditioning.
\end{rem}

 \begin{rem}\label{rem:mulambdatrival}
Let $(A,\lambda,v)\in\W$. If $(\lambda I_n-A)^*v=0$, that is, if $v$ is also a left eigenvector of $A$ with eigenvalue $\lambda$, then, from Proposition~\ref{prop:muvmulambdaexp}, one has
$$
\mu_\lambda(A,\lambda,v)=\frac{\sqrt 2}{1+\frac{|\lambda|^2}{\|A\|_F^2}}.
$$
Thus $\mu_\lambda(A,\lambda,v)\leq\sqrt{2}$. 
\end{rem}

From the previous remark we conclude that when $A$ is normal, i.e., $A^*A=AA^*$, the eigenvalue condition number $\mu_\lambda$ is not related to the distance to the discriminant variety $\Sigma$. 
On the other hand, $\mu_v$ happens to be more
interesting since, roughly speaking, $\mu_v(A,\lambda,v)$ measures how close to $\lambda$ others eigenvalues of $A$ are. More precisely, we have the following result.

\begin{lem}
Let $A$ be a normal matrix. If $(A,\lambda,v)\in\W$ then
$$
 \mu_\lambda(A,\lambda,v) \leq\sqrt{2};\qquad
\mu_v(A,\lambda,v)=\frac{\|A\|_F}{\min_{\substack{i}} |\lambda-\lambda_i| },
$$
where the minimum is taken for $\lambda_i$ eigenvalue of $A$ different from $\lambda$.
\end{lem}
\begin{proof}
The inequality for $\mu_\lambda$ follows from Remark~\ref{rem:mulambdatrival}. 
\\
Since $A$ is normal, by Proposition~\ref{prop:mulammuvunitinv}, we may assume that $A$ is the  diagonal matrix $\mbox{Diag}(\lambda,\lambda_2,\ldots,\lambda_n)$, where $\lambda,\lambda_i$ are the eigenvalues of $A$. 
Furthermore, since  $(A,\lambda,v)\in\W$, then $\lambda\neq \lambda_i$ for $i=2,\ldots n$. Thus the result follows from {Proposition~\ref{prop:muvmulambdaexp}}.
\end{proof}

\subsection{Some examples}\label{sec:examples}

In this paragraph we compute the eigenvalue and eigenvector condition numbers for some simple matrices. We denote by $e_1$
the first element of the canonical basis of the underlying $\K^n$.

(a) Let $A_1=\begin{pmatrix}1&0\\0&-1\end{pmatrix}$. Then
$(A_1,1,e_1)\in\W$. Since $A_1$ is symmetric and has eigenvalues $-1$
and $1$, we have $\mu_v(A_1,1,e_1)=1/\sqrt{2}$. 

 (b) Let $A_2=e_1 e_1^*\in\K^{n\times n}$. Then
$(A_2,1,e_1)\in\W$, where $\mu_\lambda(A_2,1,e_1)=1/\sqrt{2}$ and  $\mu_v(A_2,1,e_1)=1$. (Note that, when $n=2$, $A_2$ and $A_1+I_2$ are in the same equivalent class of
$\mathbb{P}(\K^{2\times 2})$.)

(c)  Let $B_\epsilon=\begin{pmatrix}1&0\\0&1-\epsilon\end{pmatrix}$, where $\epsilon > 0$. Then $(B_\epsilon,1,e_1)\in\W$. One has $\mu_\lambda(B_\epsilon,1,e_1)\leq\sqrt{2}$, and $\mu_v(B_\epsilon,1,e_1)={\sqrt{1+|1-\epsilon|^2}}/{\epsilon}.$ In particular, as $\epsilon\to0$ we have $\mu_v(B_\epsilon,1,e_1)\to +\infty$. 
Surprisingly,  the behaviour of $\mu_v$ can be drastically changed by the pre--conditioning procedure described in Remark~\ref{rem:muvtrans}. More precisely, the matrix 
$B_\epsilon -(\mbox{tr}(B_\epsilon)/2) I_2 =(\epsilon/2)A_1$ (where $A_1$ is given in the example (a) above), and hence $\mu_v$ attains its minimum value on the associated eigentriple.

(d)   Let 
$B_\epsilon=\begin{pmatrix}1&\epsilon\\1&1\end{pmatrix}$, where $\epsilon>0$. This matrix was
studied by Wilkinson in~\cite{Wil-ill} as an example of ill--conditioned matrix.
One has $(B_\epsilon,\lambda_\epsilon,v_\epsilon)\in\W$, where
$v_\epsilon=(\sqrt\epsilon,1)^T$ and $\lambda_\epsilon=1+\sqrt\epsilon$.
Then, we have
$\mu_v(B_\epsilon,\lambda_\epsilon,v_\epsilon)=\sqrt{3+\epsilon^2}/(2\sqrt\epsilon)$
and
$\mu_\lambda(B_\epsilon,\lambda_\epsilon,v_\epsilon)=\sqrt{1+6\epsilon+\epsilon^2}/(4\sqrt{\epsilon})$.
So, as $\epsilon$ decrease to zero, both
condition numbers growth to $+\infty$. We return to this example in Section~\ref{sec:CNT}.

(e)  Let $B_\epsilon=\begin{pmatrix}1&1/\epsilon\\0&2\end{pmatrix}$, where $\epsilon>0$. Then $(B_\epsilon,1,e_1)\in\V$. It is easily to check that $\mu_\lambda$ and $\mu_v$, at $(B_\epsilon,1,e_1)$, are larger that $1/(2\epsilon)$. Therefore both condition numbers growth to infinity as $\epsilon$ decrease to zero, even though $\chi_{B_\epsilon}(z)=(z-1)(z-2)$ is a well--posed polynomial; cf. Remark~\ref{rem:charpolyapproach}.

\subsection{Condition number of the eigenvalue problem}

The condition number of a computational problem is usually defined as the
operator norm of the map giving the first order variation of the output in terms
of the first order  variation of the input; c.f. Definition
\ref{def:eigvaleigveccn}. In our case the condition number should be the
operator norm of the condition operator  $D\mathscr{S}(A,\lambda,v)$ given in
{Section~\ref{sec:EigEigCN}}, i.e.,
$$
\|D\mathscr{S}(A,\lambda,v)\|:=\sup_{\substack{\dot B\in A^\perp\\ \|\dot B\|_F=\|A\|_F}} \|D\mathscr{S}(A,\lambda,v)\dot B\|_{(A,\lambda,v)}.
$$
However, instead of this definition, we define the condition number of the eigenvalue problem in the following way.
\begin{defn}[Condition Number of the Eigenvalue Problem]\label{def:CN}
The condition number of the eigenvalue problem is defined by
$$
\mu(A,\lambda,v):=\max\{1,\mu_v(A,\lambda,v)\},\qquad (A,\lambda,v)\in\W.
$$
\end{defn}
 In item (ii) of the next proposition we show that this definition and the
usual one are essentially equivalent.
\begin{prop}\label{prop:minmuv}
 Let $(A,\lambda,v)\in\W$. Then,
\begin{enumerate}
\item[(i)]
$
\mu_\lambda(A,\lambda,v)< ({1+\frac{|\lambda|^2}{\|A\|_F^2}})^{-1}\,
(2+\mu_v(A,\lambda,v)^2)^{1/2};
$
\item[(ii)] 
$
\mu(A,\lambda,v)< \|D\mathscr{S}(A,\lambda,v)\| <
2\, \mu(A,\lambda,v).
$
\end{enumerate}
\end{prop}
\begin{proof}
Fix a representative of $(A,\lambda,v)\in\W$ such that $\|A\|_F=1$ and $\|v\|=1$. Furthermore, by Proposition~\ref{prop:mulammuvunitinv}, without loss of generality we may assume that $v$ is the first element of the
canonical basis, and thus we may write 
$A=\begin{pmatrix}
 \lambda & a\\
0 & \hat A
\end{pmatrix}$, where $a\in\K^{1\times (n-1)}$ and $\hat
A\in\K^{(n-1)\times(n-1)}$.\\
(i) Since  $A-\lambda I_n=
\begin{pmatrix}
 0 & a\\
0 & \hat A-\lambda I_{n-1}
\end{pmatrix}$, a straightforward computation shows that $u=(1,-[(\hat A-\lambda I_{n-1})^*]^{-1} a^*)^T$ is a solution of $(A-\lambda I_n)^*u=0$, i.e., $u$ is a left eigenvector associated to $\lambda$. Here $\cdot^T$ and $\cdot^*$ denote the transpose and conjugate transpose respectively.
Then,
\begin{align*}
 \frac{|\pes vu|}{\|v\|\,\|u\|}  =
 \frac{1}{\sqrt{1+\|[(\hat A-\lambda I_{n-1})^*]^{-1}a^* \|^2}} \geq
\frac{1}{\sqrt{1+\|[(\hat A-\lambda I_{n-1})^*]^{-1}\|^2\, \|a \|^2}} 
\end{align*}
Since for every invertible matrix $B$,
$\|(B^*)^{-1}\|=\|(B^{-1})^*\|=\|B^{-1}\|$, then
$$
 \frac{|\pes vu|}{\|v\|\,\|u\|}  \geq \frac{1}{\sqrt{1+\|(\hat A-\lambda I_{n-1})^{-1}\|^2\, \|a \|^2}}.
$$
Furthermore, since
$(A,\lambda,v)\in\W$, then $|\lambda|$ and $\|\hat A\|_F$ cannot be zero at the same
time (if this is the case then $\lambda=0$ is a multiple eigenvalue). Then we have ${1=\|A\|_F=(|\lambda|^2+\|a\|^2+\|\hat A\|_F^2)^{1/2}>\|a\|}$, and therefore from Proposition
\ref{prop:muvmulambdaexp}
$$  \frac{|\pes vu|}{\|v\|\,\|u\|}>  
\frac{1}{\sqrt{1+\mu_v(A,\lambda,v)^2}}.
$$
\\
(ii): 
From the definition of the condition operator $D\mathscr{S}(A,\lambda,v)$ and
equation (\ref{eq:normadslambda}),  we obtain, for every $\dot B\in
T_{A}\pmnk$, with $\|\dot B\|_F=1$, that
$$
\| D\mathscr{S}(A,\lambda,v)(\dot B)\|_{(A,\lambda,v)}^2\geq
\frac{1}{1+\frac{|\lambda|^2}{\|A\|_F^2}} + \| D\mathscr{S}_v(A,\lambda,v)(\dot
B)\|_v^2.
$$
Then, maximizing over $\dot B\in T_{A}\pmnk$ such that $\|\dot B\|_A=1$,
we get the lower bound 
$$\| D\mathscr{S}(A,\lambda,v)\|^2\geq
\frac{1}{1+\frac{|\lambda|^2}{\|A\|_F^2}}+\mu_v(A,\lambda,v)^2.
$$
Now, the lower bound in (ii) follows from the following claim.
\\
\underline{Claim:} Let $(A,\lambda,v)\in\W$, then
$\mu(A,\lambda,v)^2<\frac{1}{1+\frac{|\lambda|^2}{\|A\|_F^2}}+\mu_v(A,\lambda,
v)^2$:\\
If $\mu_v(A,\lambda,v)\geq 1$ then
$\mu(A,\lambda,v)^2=\mu_v(A,\lambda,v)^2<\frac{1}{1+\frac{|\lambda|^2}{\|A\|_F^2
} } +\mu_v(A,\lambda,v)^2$. 
Therefore from Lemma~\ref{lem:CNlowerbound}, it suffices 
to prove the claim for the range ${\frac{1}{\sqrt{2}}\leq \mu_v(A,\lambda,v)<1}$. 
In this range ${\mu(A,\lambda,v)=1}$. Furthermore, since
$\frac{1}{1+\frac{|\lambda|^2}{\|A\|_F^2}}\geq \frac{1}{2}$ and
$\mu_v(A,\lambda,v)^2\geq \frac{1}{2}$ we reduce our problem to prove that
the last two inequalities cannot be equalities at the same time. This assertion 
follows from the fact that the condition $\|A\|_F=|\lambda|$ implies
$\mu_v(A,\lambda,v)=1$ whenever $(A,\lambda,v)\in\W$.

Let us prove the second inequality in (ii). By the definition of the condition
operator we get
$$
\| D\mathscr{S}(A,\lambda,v)\|^2
\leq \mu_v(A,\lambda,v)^2+\mu_\lambda(A,\lambda,v)^2.
$$ 
Then, from assertion (i) of this proposition we obtain that 
$$
\mu_v(A,\lambda,v)^2+\mu_\lambda(A,\lambda,v)^2< \mu(A,\lambda,v)^2+
(2+\mu_v(A,\lambda,v)^2)\leq 4\mu(A,\lambda,v)^2,
$$
proving the upper bound.
\end{proof}

\begin{rem}
Example (c) in Section~\ref{sec:examples} shows that the inequality (i) in Proposition~\ref{prop:minmuv} is far from be sharp.
\end{rem}

The next result follows immediately from Proposition~\ref{prop:mulammuvunitinv}.
\begin{prop}
The condition number $\mu$ is invariant under the action of $\mathbb{U}_n(\K)$, i.e., for every $U\in\mathbb{U}_n(\K)$, one has $\mu(UAU^{-1},\lambda,Uv)=\mu(A,\lambda,v).$
\qed
\end{prop}

\medskip

The next section is included for the sake of completeness but is not needed  for the proof of our main results.

\subsection{Condition Number Theorem}\label{sec:CNT}
In this section we study the relation of $\mu(A,\lambda,v)$ with the distance of $(A,\lambda,v)$ to
$\Sigma'$. The main objective of this section is to prove the following theorem.
\begin{thm}\label{thm:CNTmuv}
For every $(A,\lambda,v)\in\W$, 
$$
\mu_v(A,\lambda,v)\leq \frac{1}{\left(1+\frac{|\lambda|^2}{\|A\|_F^2}\right)^{1/2}}
\frac{1}{\sin(d_{\mathbb{P}^2}\left( (A,\lambda,v),\Sigma'_v)\right)},
$$
where $\Sigma'_v$ is the intersection of the fiber $\V_v=\pi_2^{-1}(v)\subset\V$ with 
the ill--posed variety $\Sigma'$.
\end{thm}
The proof of Theorem~\ref{theorem:CNT} follows immediately from Definition~\ref{def:CN} and Theorem~\ref{thm:CNTmuv}.

In general, if $(\mathbb{E},\pes{\cdot}{\cdot})$ is a finite dimensional Hermitian  vector space over $\K$, and $q:\mathbb{E}\setminus\{0\}\to\mathbb{P}(\mathbb{E})$ is the canonical quotient projection defining the projective space $\mathbb{P}(\mathbb{E})$, then, given $\Lambda\subset \mathbb{P}(\mathbb{E})$, we define  
\begin{equation}\label{def:affext}
\widehat\Lambda:=q^{-1}(\Lambda)\subset(\mathbb{E}\setminus\{0\}).
\end{equation}

Recall from the introduction that we write interchangeably a nonzero vector and its corresponding class in the projective space.

With this notation the following result is elementary.
\begin{lem}\label{lem:distap}
Given $x\in\mathbb{E}$, $x\neq 0$, and $\Lambda\subset\mathbb{P}(\mathbb{E})$ we have
$$
\sin(d_\mathbb{P}(x,\Lambda))=\frac{d_\mathbb{E}(x,\widehat\Lambda)}{\|x\|},
$$
where $d_\mathbb{E}(x,\widehat\Lambda)=\inf\{\|x-y\|:\,y\in\widehat\Lambda\}$ and $d_\mathbb{P}(x,\Lambda)=\inf\{d_\mathbb{P}(x,z):\,z\in\Lambda\}$.
\qed
\end{lem}

The next proposition is a version, adapted to this context, of a known result given by
Shub--Smale~\cite{BezIV}.

Recall that $\Sigma=\pi(\Sigma')\subset\pmnk$.

\begin{lem}\label{lem:CNT}
 Let $(A,\lambda,v)\in\W$. Then
$$
\mu_v(A,\lambda,v)=
\frac{\|A\|_F}{d_F(A-\lambda I_n,\widehat{\Sigma_{0,v}})},
$$
where $\Sigma_{0,v}=\{B\in\pmnk:\, Bv=0,\,\mbox{rank}(B^2)<n-1\}\subset\Sigma$.
\end{lem}
\begin{proof}
In Shub--Smale~\cite{BezIV} it is proved that, for a fixed triple
$(A,\lambda,v)\in\widehat\W$,
$$
d_F(\lambda I_n-A,\widehat{\Sigma_{0,v}})=\frac{1}{\|(\Pi_{v^\perp}(\lambda I_n-A)|_{v^\perp})^{-1}\|}.
$$
Then, the result follows from  {Proposition~\ref{prop:muvmulambdaexp}}.
\end{proof}

From Lemma~\ref{lem:distap} and~\ref{lem:CNT} we conclude the following result.
\begin{prop}\label{prop:cntmuv}
 Let $(A,\lambda,v)\in\W$. Then
$$
\mu_v(A,\lambda,v)= \frac{\|A\|_F}{\|A-\lambda I_n\|_F}
\frac{1}{\sin(d_{\mathbb{P}}(A-\lambda I_n,{\Sigma_{0,v}}))},
$$
where $\Sigma_{0,v}=\{B\in\pmnk:\, Bv=0,\,\mbox{rank}(B^2)<n-1\}\subset\Sigma$.
\qed
\end{prop}

\begin{rem}
 From Proposition~\ref{prop:cntmuv} and the fact that  $\sin(\cdot)\leq1$, we conclude that, if $(A,\lambda,v)\in\W$, then 
 $ \mu(A-\lambda I_n,0,v)={\sin(d_{\mathbb{P}}(A-\lambda I_n,{\Sigma_{0,v}}))}^{-1}.$
\end{rem}

\begin{proof}[Proof of Theorem~\ref{thm:CNTmuv}]
Since $(A,\lambda,v)$ and $\Sigma'_v$ are included in the fiber $\V_v$, the distance 
$d_{\mathbb{P}^2}((A,\lambda,v),\Sigma'_v)$ coincides with the projective distance of $\pi_1(A,\lambda,v)$ and $\pi_1(\Sigma'_v)$, where $\pi_1$ is the canonical projection $\pi_1:\V\to\ppuno$, that is, 
\begin{equation}\label{eq:dp2dp1CNT}
d_{\mathbb{P}^2}((A,\lambda,v),\Sigma'_v)=d_{\mathbb{P}}((A,\lambda),\pi_1(\Sigma'_v)).
\end{equation}
Note that $$
\pi_1(\Sigma'_v)=\{(B,\eta)\in\mathbb{P}(\K^{n\times n}\times \K):\; (B-\eta I_n)v=0,\;\mbox{rank}((B-\eta I_n)^2)<n-1\}.
$$

Fix a representative of $(A,\lambda,v)\in\W$, i.e., we assume $(A,\lambda,v)\in\widehat\W$. Let $d_{\K^{n\times n}\times\K}$ be the canonical distance on $\K^{n\times n}\times\K$.

\underline{Claim:}
$$
d_{\K^{n\times n}\times\K}((A,\lambda),\widehat{\pi_1(\Sigma'_v)})\leq d_F(A-\lambda I_n,\widehat{\Sigma_{0,v}}),
$$
where $\Sigma_{0,v}$ is defined in {Lemma~\ref{lem:CNT}}:\\
Since $\mathbb{U}_n(\K)$ acts by isometries on $\widehat{\V}$ (see Remark~\ref{rem:unitinv}) we may assume that $v=e_1$. Write $A=\begin{pmatrix}\lambda &a\\0&\hat A\end{pmatrix}$. Then we have
\begin{equation}\label{eq:CNTproof}
d_{\K^{n\times n}\times\K}((A,\lambda),\widehat{\pi_1(\Sigma'_v)})=\inf\{(\|A-B\|_F^2+|\lambda-\eta|^2)^{1/2}:\,(B,\eta)\in \widehat{\pi_1(\Sigma'_v)}\}. 
\end{equation}
If we consider the subset of pairs $(B,\eta)\in\widehat{\pi_1(\Sigma'_v)}$ such that $\eta=\lambda$, we get
\begin{equation*}
 \|A-B\|_F^2+|\lambda-\eta|^2 =\|(A-\lambda I_{n})-(B-\lambda I_n)\|_F^2,
\end{equation*}
where $(B-\lambda I_n)v=0$ and $\mbox{rank}((B-\lambda I_n)^2)<n-1$. Then,  
$$
d_{\K^{n\times n}\times\K}((A,\lambda),\widehat{\pi_1(\Sigma'_v)})\leq d_F(A-\lambda I_n,\widehat{\Sigma_{0,v}}),
$$
and the claim follows.

Now, from this claim and  {Lemma~\ref{lem:CNT}}, we get
$$
d_{\K^{n\times n}\times\K}((A,\lambda),\widehat{\pi_1(\Sigma'_v)})\leq\frac{\|A\|_F}{\mu_v(A,\lambda,v)}.
$$
Then, from (\ref{eq:dp2dp1CNT}) and Lemma~\ref{lem:distap}, we conclude
$$
\sin(d_{\mathbb{P}^2}((A,\lambda,v),\Sigma'_v))\leq\frac{\|A\|_F}{(\|A\|_F^2+|\lambda|^2)^{1/2}}\frac{1}{\mu_v(A,\lambda,v)}
$$
\end{proof}

\begin{rem}
When we fix a representative of $(A,\lambda,v)\in\W$, we obtain from
Proposition~\ref{prop:cntmuv} that the condition number $\mu_v(A,\lambda,v)$ is comparable to
the inverse of the sine of the projective distance of $A$ to the set of ill--posed matrices
such that $\lambda$ is not a simple eigenvalue with eigenvector $v$. However, if we remove
the last condition, the distance of $A$ to the discriminant variety $\Sigma$ could be much
smaller. This is the case of the example (d) in Section~\ref{sec:examples}. In that case, when
$\epsilon$ is small enough, $\mu_v(B_\epsilon,\lambda_\epsilon,v_\epsilon)$ has order
$\epsilon^{-1/2}$ and hence $d_{\mathbb{P}}(B_\epsilon-\lambda_\epsilon I_n,{\Sigma_{0,v}})$ has order
${\epsilon}^{1/2}$, however, the order of $d_{\mathbb{P}}(B_\epsilon,\Sigma)$ is, at least, smaller
than $\epsilon$;  cf. Wilkinson~\cite{Wil}.
\end{rem}

\medskip

\subsection{Sensitivity }

For the proof of {Theorem~\ref{teo:main}} we have to study
the rate of change of the condition number $\mu$ defined in Definition~\ref{def:CN}.

The main result of this section is the following proposition.
\begin{prop}\label{prop:mucota}
Given $\varepsilon>0$, there exist $C_\eps>0$ such that, if
$(A,\lambda,v)$, $(A',\lambda',v')$ belong to $\W$ and
$$
d_{\mathbb{P}^2}\big((A,\lambda,v),(A',\lambda',v')\big)\leq\frac{C_\eps}{\mu(A,
\lambda,v)},
$$
then
$$
\frac{\mu(A,\lambda,v)}{1+\varepsilon}\leq\mu(A',\lambda',v')\leq
(1+\varepsilon)\mu(A,\lambda,v).
$$
(One may choose
$C_\eps=\frac{\arctan\left(\frac{\varepsilon}{\sqrt{2}+\alpha
(1+\varepsilon)}\right)}{(1+\varepsilon)},$ where
$\alpha:=(1+\sqrt{5})2\sqrt{2}$.)
\end{prop}

\medskip

Before proving {Proposition
\ref{prop:mucota}} we need some additional notation.

Given $w\in\K^n$, $w\neq 0$, we define the linear operator 
\begin{equation}\label{def:extperp}
\hat{\Pi}_{w^\perp}:\mnk\to\mnk,\quad\mbox{given by}\quad
\hat{\Pi}_{w^\perp}B:=\tau_w\circ \Pi_{w^\perp}B,
\end{equation}
for every $B\in\mnk$, where  $\tau_w:w^\perp\hookrightarrow \K^n$ is the inclusion map. That is, 
$$
\hat{\Pi}_{w^\perp}Bz=Bz-\pes{Bz}{\frac{w}{\|w\|}}\frac{w}{\|w\|}.
$$

When $\mathbb{E}$ is a finite dimensional vector space
over $\K$  equipped with the Hermitian inner product $\pes\cdot\cdot$, we
define
\begin{equation}\label{def:dsubt}
d_T(w,w'):= \tan (d_\mathbb{P}(w,w')), 
\end{equation}
for all $w,\,w'\,\in\mathbb{P}(\mathbb{E})$.
 In particular,  $d_T(w,w')=\|w-w'\|_{w}$, whenever $w$ and $w'$ satisfy $\pes{w-w'}{w}=0$.

Note that $d_{\mathbb{P}}(\cdot,\cdot)\leq d_T(\cdot,\cdot)$.
Moreover, from elementary facts we have the following result.
\begin{lem}\label{lem:reldist}
Let $w,\,w'\in\mathbb{P}(\mathbb{E})$ such that
$d_\mathbb{P}(w,w')\leq \theta<\pi/2$. Then
$$
d_\mathbb{P}(w,w')\leq d_T(w,w') \leq \frac{\tan(\theta)}{\theta} \, d_\mathbb{P}(w,w')
,\;\;\mbox{for all}\; w,\,w'\in\mathbb{P}(\mathbb{E}).
\qed
$$
\end{lem}

With the notation given above we have the following result.
\begin{lem}\label{lem:deshat}
 Let $v,\,w\in\prk$ and $B\in\mnk$.
Then
$$
\left\| \hat{\Pi}_{v^\perp}B-\hat{\Pi}_{w^\perp}B \right\|\leq 2
\|B\|\, d_T(v,w).
$$
\end{lem}
\begin{proof}
Take representatives of $v$ and $w$ such that $\|v\|=1$ and
$\pes{v-w}{v}=0$. Let $u\in\K^n$, then
\begin{equation*}
\aligned
\left\|
\left(\hat{\Pi}_{v^\perp}B-\hat{\Pi}_{w^\perp}B\right)u \right\|
&= \left\|Bu -\pes{Bu}{v}v -\left(Bu -\pes{Bu}{\frac{w}{\|w\|}}\frac{w}{\|w\|}\right)   \right\| \\
&= \left\|  \pes{Bu}{\frac{w}{\|w\|}}\frac{w}{\|w\|} -\pes{Bu}{v}v  \right\| \\
&= \left\|   \pes{Bu}{ \frac{w}{\|w\|}-v }\frac{w}{\|w\|} +\pes{Bu}{v}\left(\frac{w}{\|w\|}-v \right) \right\| \\
&\leq  2 \|Bu\|\, \left\|  \frac{w}{\|w\|}-v  \right\| \leq
2 \|Bu\|\, d_T(v,w).
\endaligned
\end{equation*}
\end{proof}
\begin{notation}\label{notation:Alambda}
Given $(A,\lambda)\in\ppuno$, we denote $A_\lambda:=(\lambda I_n-A)$.
\end{notation}

\begin{rem}
Since $\left(\hat{\Pi}_{v^\perp}A_\lambda\right)v=0$ for all
$(A,\lambda,v)\in\W$, then $\left\|(\hat{\Pi}_{v^\perp}A_\lambda)^\dagger\right\|$ and 
$\left\|(\Pi_{v^\perp}A_\lambda|_{v^\perp})^{-1}\right\|$ are equal, where ${\dagger}$ denotes taking
the Moore-Penrose inverse. Then Proposition~\ref{prop:muvmulambdaexp} yields
\begin{equation}\label{eq:muvnewexp}
\mu_v(A,\lambda,v)= \|A\|_F\,
\left\|(\hat{\Pi}_{v^\perp}A_\lambda)^\dagger\right\|.
\end{equation}
\end{rem}

Let $d_{T^2}$ be the product function defined over $\pp$ by
$$
d_{T^2}((A,\lambda,v),(A',\lambda',v')):=\big(d_{T}((A,\lambda),(A',\lambda'))^2+d_{T}(v,v')^2\big)^{1/2},
$$
where $d_T$ is given in (\ref{def:dsubt}).
\begin{prop}\label{prop:cdi}
Let $\alpha:=(1+\sqrt{5})2\sqrt{2}.$  Let $(A,\lambda,v), (A',\lambda',v')\,\in\W$ such that
$$
d_{T^2}\big((A,\lambda,v),(A',\lambda',v')\big)<\frac{1}{\alpha\,\mu_v(A,\lambda,v)}.
$$
Then, the following inequality holds:
$$
\mu_v(A',\lambda',v')\leq
\frac{\left(1+\sqrt{2}d_{T^2}((A,\lambda,v),(A',\lambda',v')
)\right)\,
\mu_v(A,\lambda,v)}{1-\alpha\,\mu_v(A,\lambda,v)\,
d_{T^2}\big((A,\lambda,v),(A',\lambda',v')  \big)}.
$$
\end{prop}
\begin{proof}
Consider representatives of $(A,\lambda,v)$ and $(A',\lambda',v')$ such that: $\|A\|_F=1$, $\|v\|=1$, $(A,\lambda)-(A',\lambda')$ perpendicular to $(A,\lambda)$ in $\mnk\times\K$, and $v-v'$ perpendicular to $v$ in $\K^n$. 
From Notation~\ref{notation:Alambda} and (\ref{eq:muvnewexp}), by Wedin's Theorem (see Stewart--Sun~\cite[Theorem 3.9]{St-S}) we
have
\begin{multline*}
{\left\|
\left(\hat{\Pi}_{v^\perp}A_{\lambda}\right)^\dagger-
\left(\hat{\Pi}_{{v'}^\perp}A'_{\lambda'}\right)^\dagger
\right\|
\leq }\\
\frac{1+\sqrt 5}{2}\,
\left\|\left(\hat{\Pi}_{v^\perp}A_{\lambda}\right)^\dagger\right\|\,
\left\|\left(\hat{\Pi}_{{v'}^\perp}A'_{\lambda'}\right)^\dagger\right\| \,
\left\|\hat{\Pi}_{v^\perp}A_{\lambda}
-\hat{\Pi}_{{v'}^\perp}A'_{\lambda'}  \right\|.
\end{multline*}
Since $\left|
\left\|\left(\hat{\Pi}_{v^\perp}A_{\lambda}\right)^\dagger\right\| -
\left\|\left(\hat{\Pi}_{{v'}^\perp}A'_{\lambda'}\right)^\dagger\right\|
\right| \leq
\left\|
\left(\hat{\Pi}_{v^\perp}A_{\lambda}\right)^\dagger-\left(\hat{\Pi}_{{v'}
^\perp}A'_{\lambda'}\right)^\dagger
\right\| $,
 then,
\begin{equation*} 
\left\|
\left(\hat{\Pi}_{{v'}^\perp}A'_{\lambda'}\right)^\dagger\right\| \leq
\frac{\left\|\left(\hat{\Pi}_{v^\perp}A_{\lambda}\right)^\dagger\right\|}{1
-\frac{1+\sqrt 5}{2}\,
\left\|\left(\hat{\Pi}_{v^\perp}A_{\lambda}\right)^\dagger\right\|\,
\left\|{\hat{\Pi}_{v^\perp}A_{\lambda}-\hat{\Pi}_{{v'}^\perp}A'_{\lambda'}}
\right\| }. 
\end{equation*} 
Note that 
\begin{align}
 \left\|\hat{\Pi}_{v^\perp}A_{\lambda}-\hat{\Pi}_{{v'}^\perp}A'_{\lambda'}
\right\|  &\leq
 \left\|\hat{\Pi}_{v^\perp}A_{\lambda}-\hat{\Pi}_{{v'}^\perp}A_{\lambda}
\right\| +
\left\| \hat{\Pi}_{{v'}^\perp}A_{\lambda} -\hat{\Pi}_{{v'}^\perp}A'_{\lambda'}
\right\| \nonumber \\
&\leq   2\, \left\|A_{\lambda}\right\| \,
d_T(v,v')+\|A_{\lambda}- A'_{\lambda'} \|,\label{eq:pihatdiff}
 \end{align} 
where the second inequality follows from {Lemma~\ref{lem:deshat}} and the fact that the operator norm of $\hat{\Pi}_{v^\perp}:\mnk\to\mnk$, given in (\ref{def:extperp}), is less or equal than one. In addition,
taking into account that $(A,\lambda,v)\in\W$ and the choice of elected representatives, we get $\|A\|_F^2+|\lambda|^2\leq 2$, and  
\begin{align*}
\|A_{\lambda}- A'_{\lambda'} \| &\leq \|A-A'\| +|\lambda-\lambda'| \leq \sqrt{2}(\|A-A'\|^2 +|\lambda-\lambda'|^2)^{1/2}\\
&\leq  \sqrt{2}\,
d_T((A,\lambda),(A',\lambda'))\,\sqrt{\|A\|_F^2+|\lambda|^2}\\
&\leq  2\, d_T((A,\lambda),(A',\lambda')), 
\end{align*}
and hence from (\ref{eq:pihatdiff}), and the fact $\|A_\lambda||\leq\|A\|+|\lambda|$, we get
\begin{align*}
 \left\|\hat{\Pi}_{v^\perp}A_{\lambda}-\hat{\Pi}_{{v'}^\perp}A'_{\lambda'}
\right\| &\leq
4\, d_T(v,v')+  2\, d_T((A,\lambda),(A',\lambda'))\\
&\leq 4 ( d_T(v,v')+  d_T((A,\lambda),(A',\lambda')))\\
&\leq 4\sqrt{2}( d_T(v,v')^2+  d_T((A,\lambda),(A',\lambda'))^2)^{1/2}\\
&=4\sqrt2d_{T^2}((A,\lambda,v),(A',\lambda',v')).
\end{align*}
Then we conclude
\begin{equation}\label{eq:ineqsens}
\left\|
\left(\hat{\Pi}_{{v'}^\perp}A'_{\lambda'}\right)^\dagger\right\|
\leq
\frac{\left\|\left(\hat{\Pi}_{v^\perp}A_{\lambda}\right)^\dagger\right\|}{1
-(1+\sqrt 5)2\sqrt{2}\,
 \left\|\left(\hat{\Pi}_{v^\perp}A_{\lambda}\right)^\dagger\right\|\,
d_{T^2}((A,\lambda,v),(A',\lambda',v')) }.
\end{equation}

In addition, by the triangle inequality we have  $\|A'\|_F\leq 1+\|A-A'\|_F$. 
Then
$$
\|A-A'\|_F\leq \frac{\sqrt{2}}{(\|A\|_F^2+|\lambda|^2)^{1/2}}\|A-A'\|_F\leq \sqrt{2} d_T((A,\lambda),(A',\lambda')),
$$
and hence $\|A'\|_F\leq
1+\sqrt{2}d_T((A,\lambda),(A',\lambda'))$. Then, the proof follows by multiplying both sides of equation~(\ref{eq:ineqsens}) by $\|A'\|_F$.
\end{proof}

\begin{prop}\label{prop::musensdt2}
 Given $\varepsilon>0$, there exist $c_\varepsilon>0$ such that, if
$(A,\lambda,v)$, $(A',\lambda',v')$ lie in $\W$ and
$$
d_{T^2}\big((A,\lambda,v),(A',\lambda',v')\big)\leq\frac{c_\varepsilon}{\mu(A,
\lambda,v)},
$$
then
$$
\mu(A',\lambda',v')\leq (1+\varepsilon)\mu(A,\lambda,v).
$$
(One may choose
$\displaystyle{c_\varepsilon=\frac{\varepsilon}{\sqrt{2}+\alpha
(1+\varepsilon)}}$, where $\alpha=(1+\sqrt 5)2\sqrt2$.)
\end{prop}
\begin{proof}
The condition
$$
d_{T^2}\big((A,\lambda,v),(A',\lambda',v')\big)\leq\frac{c}{\mu(A,\lambda,v)},
$$
implies
$$
d_{T^2}\big((A,\lambda,v),(A',\lambda',v')\big)\leq\frac{c}{\mu_v(A,\lambda,v)}
.
$$
From {Proposition~\ref{prop:cdi}} and the fact that $\mu\geq1$, if 
$c<1/\alpha$ and
$$
\frac{1+\sqrt{2} c}{1-\alpha c}\leq 1+\varepsilon,
$$
we get 
$$
\mu_v(A',\lambda',v')\leq (1+\varepsilon)\mu_v(A,\lambda,v).
$$
Then 
\begin{align*}
\mu(A',\lambda',v')&=\max\{1,\mu_v(A',\lambda',v')\} \\
&\leq \max\{1,(1+\varepsilon)\mu_v(A,\lambda,v)\}\\
&\leq 
(1+\varepsilon)\max\{1,\mu_v(A,\lambda,v)\}=(1+\varepsilon)\mu(A,
\lambda,v).
\end{align*}

One may choose
$\displaystyle{c_\varepsilon=\frac{\varepsilon}{\sqrt{2}+\alpha
(1+\varepsilon)}}$.

\end{proof}

\begin{cor}\label{cor:mucota0}
 Given $\varepsilon>0$, there exist $c_\varepsilon'>0$ such that, if $(A,\lambda,v)$, $(A',\lambda',v')$ lie in $\W$ and
$$
d_{\mathbb{P}^2}\big((A,\lambda,v),(A',\lambda',v')\big)\leq\frac{
c_\varepsilon'}{\mu(A,\lambda,v)},
$$
then
$$
\mu(A',\lambda',v')\leq (1+\varepsilon)\mu(A,\lambda,v).
$$
(One may choose
$c_\varepsilon'=\arctan\left(\frac{\varepsilon}{\sqrt{2}+\alpha
(1+\varepsilon)}\right)$  where
$\alpha:=(1+\sqrt{5})2\sqrt{2}$.)
\end{cor}
\begin{proof}
 By Lemma~\ref{lem:reldist}, if
$$
d_{\mathbb{P}^2}\big((A,\lambda,v),(A',\lambda',v')\big)\leq\frac{c'}{\mu(A,
\lambda,v)},
$$
then
\beqna
d_{T^2}\big((A,\lambda,v),(A',\lambda',v')\big) & \leq &
\frac{\tan(c')}{c'}\,
d_{\mathbb{P}^2}\big((A,\lambda,v),(A',\lambda',v')\big)\\
&\leq&
\frac{\tan(c')}{\mu(A,\lambda,v)},
\eeqna
so we just need to choose $c'$ such that $\tan(c')\leq c_\varepsilon$ from
Proposition~\ref{prop::musensdt2}.
\end{proof}

\begin{proof}[Proof of Proposition~\ref{prop:mucota}]
From {Corollary~\ref{cor:mucota0}}, there exist $c'>0$ such
that, if $(A,\lambda,v)$, $(A',\lambda',v')$ $\in\W$ are such that
$$d_{\mathbb{P}^2}\big((A,\lambda,v),(A',\lambda',v')\big)\,\mu(A,\lambda,
v)\leq c',
$$
then
$$
\mu(A',\lambda',v')\leq(1+\varepsilon)\mu(A,\lambda,v).
$$
It is enough to take $c'$ such that $c'\leq
\arctan\left(\frac{\varepsilon}{\sqrt{2}+\alpha
(1+\varepsilon)}\right)$.
In this case we have
$$
d_{\mathbb{P}^2}\big((A,\lambda,v),(A',\lambda',v')\big)\,\mu(A',\lambda',
v')\leq c'(1+\varepsilon).$$
Then, by the same argument,  if
$c'(1+\varepsilon)\leq\arctan\left(\frac{\varepsilon}{\sqrt{2}+\alpha
(1+\varepsilon)}\right)$ we have the other inequality.
\end{proof}
\medskip

\section{Newton's Method}\label{sec:NM}
\subsection{Introduction}
In this section we start describing the Newton method defined in Section~\ref{sec:BHNMeth}.
The main goal of this section is to prove  {Theorem~\ref{teo:newton}}.

Let us recall some definitions from the introduction.

Given a nonzero matrix $A\in\mnk$, let $F_A:\K\times\K^n \to \K^n$ be the
evaluation map
$$
F_A(\lambda,v):=(\lambda I_n-A)v.
$$

This map is homogeneous of degree $1$ in $v$. 
Its derivative ${DF_A(\lambda,v):\K\times\K^n\to\K^n}$ satisfies
\begin{equation}\label{eq:DFAalfa}
DF_A(\lambda,\alpha v)(\dot\lambda,\alpha \dot v)=\alpha DF_A(\lambda,
v)(\dot\lambda,\dot v),
\end{equation}
for all $(\dot\lambda,\dot v)\in\K\times\K^n$, and nonzero scalar $\alpha$.

\begin{defn}
 Given a nonzero matrix $A\in\mnk$, we define the \textit{Newton map associated
to} $A$ to be the map $N_{A}:\K\times(\K^n\setminus\{0\})\to\K\times(\K^n\setminus\{0\})$ given by 
$$
N_{A}(\lambda,v):=(\lambda,v)-\big(DF_A(\lambda,v)|_{\K\times
v^{\perp}}\big)^{-1}F_A(\lambda,v),
$$
defined for all $(\lambda,v)$ such that $DF_A(\lambda,v)|_{\K\times
v^{\perp}}$ is invertible.
\end{defn}
Note that, from (\ref{eq:DFAalfa}), the map $N_A$  induces a map from $\K\times\prk$ into itself (defined almost everywhere); cf. Remark~\ref{rem:Nwelldefpp}.

\begin{lem}\label{lem:DFAinv}
Let  $A\in\mnk$ be a nonzero matrix and $(\lambda,v)\in\K\times\K^n$, $v\neq0$. 
The map $N_A$ is well--defined at $(\lambda,v)$ if and only if $\Pi_{v^\perp}(\lambda I_n-A)|_{v^\perp}$ is invertible.
\end{lem}
\begin{proof}
The map $N_A$ is well--defined provided that the linear operator $DF_A(\lambda,v)|_{\K\times
v^\perp}$, from $\K\times v^\perp$ into $\K^n$, is invertible. Differentiating $F_A$ with respect to $\lambda$ and $v$ yields
$$ 
 DF_A(\lambda,v)(\dot\lambda,\dot v)=\dot\lambda v+(\lambda I_n-A)\dot v, \quad (\dot\lambda,\dot v)\in\K\times\K^n.
 $$ 
Fix a basis of $\K^n$ and let $w\in\K^n$. 
Solving the linear equation 
$ DF_A(\lambda,v)(\dot\lambda,\dot v)=w$, 
for $(\dot\lambda,\dot v)\in\K\times v^\perp$, is equivalent to solve the system of equations:
\begin{equation}\label{eq:matexpNA}
 \begin{pmatrix}v & \lambda I_n-A\\ 0& v^*\end{pmatrix} \begin{pmatrix}\dot\lambda\\ \dot v\end{pmatrix}=\begin{pmatrix}w\\ 0\end{pmatrix}, \qquad\mbox{for}\quad (\dot\lambda,\dot v)\in\K\times\K^n.
\end{equation}
Hence $N_A(\lambda,v)$ is well defined if and only if the matrix given in (\ref{eq:matexpNA}) is invertible.

Let $U\in\mathbb{U}_n(\K)$ such that $Uv=\|v\|e_1$. Then,
$$
\begin{pmatrix}U & 0\\ 0& 1\end{pmatrix} \begin{pmatrix}v & \lambda I_n-A\\ 0& v^*\end{pmatrix} \begin{pmatrix}1 & 0\\ 0& U^*\end{pmatrix} =\begin{pmatrix}\|v\|e_1 & U(\lambda I_n-A)U^*\\ 0& \|v\|e_1^*\end{pmatrix}.
$$
Now,  expanding the determinant of $\begin{pmatrix}\|v\|e_1 & U(\lambda I_n-A)U^*\\ 0& \|v\|e_1^*\end{pmatrix}$ by the first column, and thereafter by the last row, we conclude that the matrix $ \begin{pmatrix}v & \lambda I_n-A\\ 0& v^*\end{pmatrix}$ is invertible if and only if the operator $\Pi_{v^\perp}(\lambda I_n-A)|_{v^\perp}$ is invertible.
\end{proof}

\begin{rem}\label{rem:expNA}
From the proof of Lemma~\ref{lem:DFAinv} we obtain that $N_A$ has the simple matrix expression 
$$
N_A \begin{pmatrix}\lambda\\ v\end{pmatrix}=
 \begin{pmatrix}\lambda\\ v\end{pmatrix}-
 \begin{pmatrix}v & \lambda I_n-A\\ 0& v^*\end{pmatrix}^{-1}
 \begin{pmatrix}(\lambda I_n-A)v\\ 0\end{pmatrix}.
$$ 
Furthermore, solving the system (\ref{eq:matexpNA}) for $w=(\lambda I_n-A)v$,  we
conclude that if  $\Pi_{v^\perp}(\lambda I_n-A)|_{v^\perp}$ is invertible then
the map $N_A$ is given by $N_A(\lambda,v)=(\lambda-\dot\lambda,v-\dot v),$ where
$$
\dot v  =\left(\Pi_{v^\perp}(\lambda
I_n-A)\big|_{v^\perp}\right)^{-1}\Pi_{v^\perp}(\lambda I_n-A)v; \qquad
\dot\lambda =\frac{\pes{(\lambda I_n-A)(v-\dot v)}{v}}{\pes vv}.
$$
\end{rem}

\begin{defn}\label{defn:newseq-approxzero}
Let $A\in\mnk$ be a nonzero matrix, and let $(\lambda_0,v_0)$ in ${\K\times\prk}$. 
 We say that the triple $(A,\lambda_0,v_0)$ is an 
\emph{approximate solution} of the eigenvalue problem ${(A,\lambda,v)\in \V}$, if the sequence 
$(A,N_A^k(\lambda_0,v_0))$, $k=0,1,\ldots$ is defined and satisfies
$$
 d_{\mathbb{P}^2}\left((A,N_A^k(\lambda_0,v_0)),(A,\lambda,v)\right)\leq 
 \left(\frac{1}{2}\right)^{2^k-1} d_{\mathbb{P}^2}\left((A,\lambda_0,v_0),(A,\lambda,v)\right),
 $$
 for all positive integers $k$.
\end{defn}
Recall from Remark~\ref{rem:AZwelldefpp} that the notion of approximate solution, and the  sequence $(A,N_A^k(\lambda_0,v_0))$, $k=0,1,\ldots$, are well--defined on $\pp$.

\medskip

\subsection{Approximate Solution Theorem}
The main tool to prove Theorem~\ref{teo:newton} is the following result.
\begin{prop}\label{prop:NewtonAfin}
Let $0<c\leq 1/(2\sqrt{2})$. Let $A\in\mnk$ such
that $\|A\|_F=1$, and let $(\lambda,v),\,(\lambda_0,v_0)\in\K\times\prk$. 
If $(A,\lambda,v)\in\W$ and
$$(|\lambda_0-\lambda|^2+d_\mathbb{P}(v_0,v)^2)^{1/2}<\frac{c}{ \mu(A,\lambda,v)},$$
 then, the sequence $(\lambda_k,v_k):=N^k_A(\lambda_0,v_0)$ satisfies
$$
(|\lambda_k-\lambda|^2+d_\mathbb{P}(v_k,v)^2)^{1/2}\leq
 \left(\frac{2\tan(c)}{1-\sqrt{2}\,c}\right)\,
  \left(\frac{1}{2}\right)^{2^k-1}\,
(|\lambda_0-\lambda|^2+d_\mathbb{P}(v_0,v)^2)^{1/2},
$$
for all positive integers $k$.
\end{prop}
(Since we do not find an appropriate version in the literature to cite, we include a proof of this proposition in the appendix of this paper.)

\begin{rem}
Some expressions given in Proposition~\ref{prop:NewtonAfin} are not scale invariant in $(A,\lambda)$ (and $(A,\lambda_0)$), and thus a restriction on $\|A\|$ is required.
\end{rem}

Picking $c$ in Proposition~\ref{prop:NewtonAfin} such that $0<c\leq 1/(2\sqrt{2})$ and
${2\tan(c)/(1-\sqrt{2}c)\leq 1}$  we have the following result, which is interesting \textit{per se}.
\begin{thm}\label{teo:gammadp}
There is a universal constant $c_0>0$ with the following property.
Let $A\in\mnk$ such
that $\|A\|_F=1$, and let $(\lambda,v),\,(\lambda_0,v_0)\in\K\times\prk$. 
If $(A,\lambda,v)\in\W$ and 
 $$(|\lambda_0-\lambda|^2+d_\mathbb{P}(v_0,v)^2)^{1/2}<\frac{c_0}{\mu(A,\lambda,v)},$$
  then the sequence $(\lambda_k,v_k):=N_A^k(\lambda_0,v_0)$ satisfies
$$
(|\lambda_k-\lambda|^2+d_\mathbb{P}(v_k,v)^2)^{1/2}\leq
  \left(\frac{1}{2}\right)^{2^k-1}\,
(|\lambda_0-\lambda|^2+d_\mathbb{P}(v_0,v)^2)^{1/2},
$$
for all positive integers $k$. (One may choose $c_0=0.288$.)\qed
\end{thm}

This theorem is a version,  for the map $N_A:\K\times\prk\to\K\times\prk$,  of a well--known
theorem in the literature, namely, the \emph{Smale $\gamma$-Theorem} (or \emph{Approximate Solution Theorem}), which gives  the size of
the basin of attraction of Newton's method; see Blum et al.~\cite[Theorem 1, pp. 263]{B-C-S-S}.

\subsection{Proof of Theorem~\ref{teo:newton}}

For the proof of Theorem~\ref{teo:newton} we need a technical lemma. Its
proof is included in the appendix.

\begin{lem}\label{des:dt}
Let $A\in\mnk$ such
that $\|A\|_F=1$, and let $(\lambda,v),\,(\lambda',v')\in\K\times\prk$. 
\begin{enumerate}
\item  If $|\lambda-\lambda'|\leq c<\sqrt 2$, then,
$$
 d_{\mathbb{P}^2}\big((A,\lambda,v),(A,\lambda',v')\big)
 \leq \beta_c\,\left(|\lambda-\lambda'|^2+ d_{\mathbb{P}}(v,v')^2\right)^{1/2},
$$
where $\beta_c=(1-c^2/2)^{-1/2}$.
\item If
$d_{\mathbb{P}^2}\big((A,\lambda,v),(A,\lambda',v')\big)
<\theta<\pi/4$, then,
$$
(|\lambda-\lambda'|^2+d_T(v,v')^2)^{1/2}\leq R_\theta\,
d_{\mathbb{P}^2}((A,\lambda.v),(A,\lambda',v')),
$$
where $R_\theta=[\sqrt{2}/\cos(\theta+\pi/4)^3]^{1/2}$.
\end{enumerate}
\end{lem}

Let $\theta_0$ such that $R_{\theta_0}\,\theta_0=1/(2\sqrt{2})$, where
$R_\theta$ is given in {Lemma~\ref{des:dt}} (${\theta_0\approx 0.1389}$).
\begin{prop}\label{prop:AZTA}
Let $0<c\leq\theta_0$. 
Let $A\in\mnk$ be a nonzero matrix, and let $(\lambda,v),\,(\lambda_0,v_0)\in\K\times\prk$. 
If $(A,\lambda,v)\in\W$ and
$$d_{\mathbb{P}^2}\big((A,\lambda_0,v_0),(A,\lambda,v)\big)<\frac{c}{\mu(A,\lambda,v)},$$
 then the sequence $(A,N_A^k(\lambda_0,v_0))$, $k=0,1,\ldots$ satisfies
\begin{multline*}
 d_{\mathbb{P}^2}\big((A,N_A^k(\lambda_0,v_0)),(A,\lambda,v)\big)
 \leq  \\
 \qquad \quad \leq R_{c}\,\beta_{c\, R_c}\left(\frac{2\tan(c\,R_c)}{1-\sqrt{2}\,c\,R_c}\right)
\,\left(\frac{1}{2}\right)^{2^k-1}
d_{\mathbb{P}^2}\big((A,\lambda_0,v_0),(A,\lambda,v)\big),
\end{multline*}
for all positive integers $k$, where $\delta(c):= c/(1-c)$.
\end{prop}
\begin{proof}
From Remark~\ref{rem:Nwelldefpp}, we may assume $\|A\|_F=1$.

Since $d_\mathbb{P}(\cdot)\leq d_T(\cdot)$, by Lemma~\ref{des:dt}--(2), one has
\beq
(|\lambda_0-\lambda|^2+d_\mathbb{P}(v_0,v)^2)^{1/2}\leq \frac{c
R_c}{\mu(A,\lambda,v)}.
 \eeq
Since $c\leq\theta_0$, we have $c
R_c\leq 1/(2\sqrt{2})$, and then Proposition~\ref{prop:NewtonAfin} yields
\begin{multline}
(|\lambda_k-\lambda|^2+d_\mathbb{P}(v_k,v)^2)^{1/2}\leq \label{eq:distnewtaf}\\
 \left(\frac{2\tan(c\,R_c)}{1-\sqrt{2}\,c\,R_c}\right)\,\left(\frac{1}{2}\right)^{2^k-1}\, (|\lambda_0-\lambda|^2+d_\mathbb{P}(v_0,v)^2)^{1/2},
 \end{multline}
for all $k>0$, where $(\lambda_k,v_k):=N^k_A(\lambda_0,v_0)$.
Since ${(|\lambda_0-\lambda|^2+d_\mathbb{P}(v_0,v)^2)^{1/2} < c\, R_c}$, we
deduce from Lemma~\ref{des:dt} and (\ref{eq:distnewtaf}) that
\begin{align*}
& d_{\mathbb{P}^2}((A,\lambda_k,v_k),(A,\lambda,v))\leq \\
&\qquad\qquad \leq \beta_{c\, R_c}
\left(\frac{2\tan(c\,R_c)}{1-\sqrt{2}\,c\,R_c}\right)\,\left(\frac{1}{2}
\right)^{2^k-1}\,
(|\lambda_0-\lambda|^2+d_\mathbb{P}(v_0,v)^2)^{1/2}
\\
&\qquad\qquad \leq R_{c}\,\beta_{c\, R_c}
\left(\frac{2\tan(c\,R_c)}{1-\sqrt{2}\,c\,R_c}\right)\,\left(\frac{1}{2}
\right)^{2^k-1}\, d_{\mathbb{P}^2}((A,\lambda_0,v_0),(A,\lambda,v)).
\end{align*}
 (Note that $c\leq \theta_0<\frac{\pi}{4}$.)
\end{proof}

\subsubsection{Proof of Theorem~\ref{teo:newton}}
\begin{proof}[Proof of {Theorem~\ref{teo:newton}}]
From {Proposition~\ref{prop:AZTA}}, proof of {Theorem~\ref{teo:newton}} follows picking ${c_0>0}$ such that $c_0\leq\theta_0$ and
$R_{c_0}\,\beta_{c_0\, R_{c_0}} \left(\frac{2\tan(c_0\,R_{c_0})}{1-\sqrt{2}\,c_0\,R_{c_0}}\right)\leq 1$. (One may choose $c_0=0.0739$.)
\end{proof}

\medskip

\section{Complexity bound}\label{sec:proofmain}

\subsection{Condition length}\label{sec:condlengh}

Let us start recalling some basic definition. 

Let $\mathbb{E}$ be a finite dimensional Hilbert space. A function $\alpha:[0,1]\to\mathbb{E}$ is an \textit{absolutely continuous path} if it is almost everywhere differentiable,  its derivative $\dot\alpha(t)$ is an integrable function, and
$$
\alpha(t)=\alpha(0)+\int_0^t\dot\alpha(s)\,ds.
$$
We say that the (projective) path $\alpha:[0,1]\to\mathbb{P}(\mathbb{E})$ is an absolutely continuous path if it is the projection, under the quotient canonical map $\mathbb{E}\setminus\{0\}\to \mathbb{P}(\mathbb{E})$, of an absolutely continuous path in $\mathbb{E}\setminus\{0\}$. 

Let us recall some definition from the introduction.

\begin{defn}
The \textit{condition length} of an absolutely continuous path ${\Gamma:[0,1]\to\W}$ is defined by
$$
\ell_\mu(\Gamma):=\int_a^b \|\dot \Gamma(t)\|_{\Gamma(t)}\,\mu(\Gamma(t))\,dt.
$$ 
\end{defn}

The next proposition is useful for the proof of our main theorem.

\begin{prop}\label{prop:cotak}
Given $\varepsilon>0$, $C_\eps>0$ as in {Proposition~\ref{prop:mucota}}, and
$\Gamma:[0,1]\to\W$ an absolutely continuous path (with
$\ell_\mu(\Gamma)<\infty$), define the sequence 
$t_0,\,t_1,\ldots$
in
$[0,1]$ such that:\\
$\bullet$ $t_0=0$;
\\
$\bullet$ $t_k$ such that
$\mu(\Gamma(t_{k-1}))\int_{t_{k-1}}^{t_k}\|\dot\Gamma(s)\|_{\Gamma(s)}
ds=C_\eps$,\\
 whenever
$\mu(\Gamma(t_{k-1}))\int_{t_{k-1}}^{1}\|\dot\Gamma(s)\|_{\Gamma(s)}ds>C_\eps$;
\\
$\bullet$ else define $t_k=t_K=1$.
\\
Then,
$$
K\leq \frac{1+\varepsilon}{C_\eps}\, \ell_\mu(\Gamma)+1.
$$
\end{prop}
\begin{proof}
Whenever $k<K$ (where $K\in\N\cup\{\infty\}$), given $t\in [t_{k-1},t_k]$,
$$
d_{\mathbb{P}^2}(\Gamma(t_{k-1}),\Gamma(t)) \leq
\int_{t_{k-1}}^{t_k}\|\dot\Gamma(s)\|_{\Gamma(s)}ds =
\frac{C_\eps}{\mu(\Gamma(t_{k-1}))}.
$$
By the first inequality in Proposition~\ref{prop:mucota}, we get
$$
\int_{t_{k-1}}^{t_k}\|\dot\Gamma(s)\|_{\Gamma(s)}\mu(\Gamma(s))\,ds\geq \frac{\mu(\Gamma(t_{k-1}))}{1+\varepsilon}\,\int_{t_{k-1}}^{t_k}\|\dot\Gamma(s)\|_{\Gamma(s)}dt=\frac{C_\eps}{1+\varepsilon}.
$$
Since $\ell_\mu(\Gamma)<\infty$, then $K<\infty$, and adding, yields
$$
\ell_\mu(\Gamma)\geq (K-1)\frac{C_\eps}{1+\varepsilon}.
$$
\end{proof}

\subsection{Proof of Theorem~\ref{teo:main}}

\begin{proof}[{Proof of Theorem~\ref{teo:main}:}]

Since $\Gamma(t)=(A(t),\lambda(t),v(t))$, $0\leq t\leq 1$, is an absolutely continuous path in $\W$, we may assume that the path $A(t)$, in $\mnk\setminus\{0_n\}$, is absolutely continuous. In addition, without loss of generality we may assume that $\|A(t)\|_F=1$, for every $t\in[0,1]$.

The idea of the proof is to show that the mesh $0=t_0<t_1<\cdots < t_K=1$ given in Proposition~\ref{prop:cotak}, for some $\varepsilon>0$ to be defined afterwards, guarantees that the predictor--corrector sequence $\hat{\Gamma}(t_{k}):=(A(t_{k}),\lambda_{k},v_{k})$, where 
 $$
 (\lambda_{k+1},v_{k+1}):=N_{A(t_{k+1})}(\lambda_k,v_k),\quad 0\leq k\leq K-1,
 $$
 approximates the path $\Gamma$ provided that $(A(0),\lambda_0,v_0)$ is an approximate solution of $\Gamma(0)$.

The proof is by induction.

Let $\varepsilon>0$, and let $C_\varepsilon$ as in Proposition~\ref{prop:mucota}.
Assume that $\Gamma(t_k),\,\hat{\Gamma}(t_{k}),\,\Gamma(t_{k+1})$ are such that,
$$
d_{\mathbb{P}^2}({\Gamma}(t_{k}),{\Gamma}(t_{k+1}))\leq 
\frac{C_\varepsilon}{\mu({\Gamma}(t_{k}))},\quad\mbox{and}\quad
d_{\mathbb{P}^2}({\Gamma}(t_{k}),\hat{\Gamma}(t_{k}))\leq  
\frac{C_\varepsilon}{\mu({\Gamma}(t_{k}))}.
$$
Then, 
\begin{align*}
&d_{\mathbb{P}^2}({\Gamma}(t_{k+1}),(A(t_{k+1}),\lambda_k,v_k)) 
 \leq \\
&\qquad\quad \leq d_{\mathbb{P}^2}({\Gamma}(t_{k+1}),{\Gamma}(t_{k}))+ 
d_{\mathbb{P}^2}({\Gamma}(t_{k}),\hat{\Gamma}(t_{k}))+ 
d_{\mathbb{P}^2}(\hat{\Gamma}(t_{k}),(A(t_{k+1}),\lambda_k,v_k))
\\
&\qquad\quad<
\frac{2C_\varepsilon}{\mu({\Gamma}(t_{k}))} +
d_{\mathbb{P}^2}(\hat{\Gamma}(t_{k}),(A(t_{k+1}),\lambda_k,v_k)).
\end{align*}
Note that
$$
d_{\mathbb{P}^2}(\hat{\Gamma}(t_{k}),(A(t_{k+1}),\lambda_k,v_k))= d_{\mathbb{P}}((A(t_{k}),\lambda_k),(A(t_{k+1}),\lambda_k)).
$$

\underline{Claim:} One has
$$
d_{\mathbb{P}}((A(t_{k}),\lambda_k),(A(t_{k+1}),\lambda_k))
\leq d_{\mathbb{P}}(A(t_{k}),A(t_{k+1})):
$$
For the ease of notation  let us denote $a:=A(t_k)$, ${a':=A(t_{k+1})}$, ${\lambda:=\lambda_k}$,
 ${\theta_0:=d_{\mathbb{P}}(a,a')}$, and 
$\theta_\lambda:=d_{\mathbb{P}}
((a,\lambda),(a',\lambda))$. Since $\|A(t)\|_F=1,\,
0\leq t\leq 1$, we have ${\|a\|_F=\|a'\|_F=1}$. 
In addition, by the law of cosines we have that 
$\cos\theta_0=1-\|a-a'\|_F^2/2$, and $\cos\theta_\lambda=1-\|a-a'\|_F^2/(2(1+|\lambda|^2)).$
Then
$\cos\theta_0\leq\cos\theta_\lambda$. Since 
$\theta_0,\theta_\lambda\,\in[0,\pi]$,
we conclude $\theta_\lambda\leq\theta_0$.

Furthermore,
\begin{align*}
 d_{\mathbb{P}}(A(t_{k}),A(t_{k+1}))
&\leq \int_{t_k}^{t_{k+1}} \|\dot{A}(s) \|_{A(s)}\,ds \\
&\leq \sqrt{2}\int_{t_k}^{t_{k+1}}
\|D\mathscr{S}_{\lambda}(\Gamma(t))\dot{A}(s) \|_{(A(s),\lambda(s))}\,ds \\
&\leq \sqrt{2}\int_{t_k}^{t_{k+1}} \|\dot\Gamma(s)\|_{\Gamma(s)}\,ds,
\end{align*}
where the second inequality follows from the trivial lower bound which one may
obtain from (\ref{eq:normadslambda}) and the assumption $\|A(s)\|_F=1$ (and
hence $\pes{\dot A(s)}{A(s)}_F=0$).

Since, by construction, $\int_{t_k}^{t_{k+1}} \|\dot\Gamma(s)\|_{\Gamma(s)}\,ds\leq C_\varepsilon/\mu({\Gamma}(t_{k}))$  we conclude
 $$
d_{\mathbb{P}^2}({\Gamma}(t_{k+1}),(A(t_{k+1}),\lambda_k,v_k)) <
\frac{(2+\sqrt{2})C_\varepsilon}{\mu({\Gamma}(t_{k}))}.
$$
Furthermore, since $d_{\mathbb{P}^2}({\Gamma}(t_{k}),{\Gamma}(t_{k+1}))< 
{C_\varepsilon}{\mu({\Gamma}(t_{k}))}^{-1}$, Proposition
\ref{prop:mucota} yields
$$
d_{\mathbb{P}^2}({\Gamma}(t_{k+1}),(A(t_{k+1}),\lambda_k,v_k)) \leq
\frac{(1+\varepsilon)(2+\sqrt{2})C_\varepsilon}{\mu({\Gamma}(t_{k+1}))}.
$$
From {Proposition~\ref{prop:AZTA}}, if
$c:=(1+\varepsilon)C_\varepsilon(2+\sqrt{2})\leq\theta_0$, then
\begin{align*}
\lefteqn{d_{\mathbb{P}^2}(\hat\Gamma(t_{k+1}),{\Gamma}(t_{k+1}))
\leq} \\
&\qquad\qquad  \leq R_{c}\,\beta_{c\, R_c} \left(\frac{2\tan(c\,R_c)}{1-\sqrt{2}\,c\,R_c}\right)\frac{1}{2}\,  d_{\mathbb{P}^2}((A(t_{k+1}),\lambda_k,v_k),\Gamma(t_{k+1})) \\
&\qquad\qquad \leq \frac{R_{c}\,\beta_{c\, R_c} \left(\frac{2\tan(c\,R_c)}{1-\sqrt{2}\,c\,R_c}\right)\frac{1}{2}\,c }{\mu({\Gamma}(t_{k+1}))}.
 \end{align*}
Then, if $\varepsilon$ is small enough such that $c\leq\theta_0$ and $R_{c}\,\beta_{c\, R_c} \left(\frac{2\tan(c\,R_c)}{1-\sqrt{2}\,c\,R_c}\right)\frac{1}{2}\,c <
C_\varepsilon$, we get
$$
d_{\mathbb{P}^2}(\Gamma(t_{k+1}),\hat{\Gamma}(t_{k+1}))\leq 
\frac{C_\varepsilon}{\mu({\Gamma}(t_{k}))}.
$$
Moreover, if the $\varepsilon$ picked above also satisfies $C_\varepsilon\leq c_0$ (where $c_0$ is given in Theorem~\ref{teo:newton}), then we have concluded that $\hat{\Gamma}(t_{k+1})$ is an
approximate solution of $\Gamma(t_{k+1})$ provided that $\hat{\Gamma}(t_{k})$ is an
approximate solution of $\Gamma(t_{k})$. This just finishes the induction step.
 (One can choose $\eps=0.1640$, $C_\eps\approx 0.01167$, and $C=100$.)
\end{proof}

\bigskip

\section{Appendix}
This section is divided in two parts.
In the first part we include a proof of {Proposition~\ref{prop:NewtonAfin}}. In the second,  we prove {Lemma~\ref{des:dt}}.

\subsection{Proof of Proposition~\ref{prop:NewtonAfin}}

Given a nonzero matrix $A\in\mnk$, recall that  evaluation map $F_A:\K\times\K^n \to \K^n$ is given by 
$F_A(\lambda,v):=(\lambda I_n-A)v.$

Throughout this section, we consider the canonical Hermitian structure on
${\K\times\K^n}$.

\subsubsection{Preliminaries and technical lemmas}
The next result follows by elementary computations.
\begin{lem}\label{lem:tecvvpr}
Let $v,\,v'\in\prk$ such that $d_{\mathbb{P}}(v,v')<\pi/2$. Let ${\Pi_{v^{\perp}}|}_{{v'}^{\perp}}:{v'}^{\perp}\to
v^{\perp}$ be the restriction of the orthogonal projection
$\Pi_{v^{\perp}}$ of $\K^n$ onto ${v'}^{\perp}$. Then,
$$
\|\left({\Pi_{v^{\perp}}|}_{{v'}^{\perp}}\right)^{-1}\|=\frac{1}{\cos (d_{\mathbb{P}}(v,v'))}.
\qed
$$
\end{lem}
(In  the preceding lemma, we consider the spaces $v^\perp$ and 
$v'^\perp$ as subspaces of $\K^n$ with the canonical Hermitian
structure.)

\begin{lem}\label{lem:tec}
Let $(A,\lambda,v)\in\W$ and $v'\in\prk$ such that $d_{\mathbb{P}}(v,v')<\pi/2$.
\begin{enumerate}
\item[(i)] For every $(\dot\lambda,\dot v)\in\K\times v^\perp$ we have
$$
\left(DF_A(\lambda,v)|_{\K\times {v'}^{\perp}}\right)^{-1}\, DF_A(\lambda,v)|_{\K\times v^{\perp}}(\dot\lambda,\dot v)=
(\dot\lambda, ({\Pi_{v^{\perp}}|}_{{v'}^{\perp}})^{-1}(\dot v)).
$$
\item[(ii)]  
$$
\|\left(DF_A(\lambda,v)|_{\K\times {v'}^{\perp}}\right)^{-1}\, DF_A(\lambda,v)|_{\K\times v^{\perp}}\|=
 \frac{1}{\cos (d_{\mathbb{P}}(v,v'))}; 
 $$
\item[(iii)]
$$
\|\left(DF_A(\lambda,v)|_{\K\times {v'}^{\perp}}\right)^{-1}\|\leq \frac{\|\left(DF_A(\lambda,v)|_{\K\times v^{\perp}}\right)^{-1}\|}{\cos(d_{\mathbb{P}}(v,v'))}.
$$
\end{enumerate}
\end{lem}
\begin{proof}
(i): 
Given $(\dot\lambda,\dot v)\in \K\times v^{\perp}$, let $(\dot\eta,\dot w)\in\K\times {v'}^\perp$ such that
$$
\big(\dot\eta,\dot w\big) = \left(DF_A(\lambda,v)|_{\K\times {v'}^{\perp}}\right)^{-1} \, DF_A(\lambda,v)|_{\K\times v^{\perp}}(\dot\lambda,\dot v).
$$
Then,
$$
\dot\eta v+ (\lambda I_n-A)\dot w= \dot \lambda v + (\lambda I_n-A)\dot v.
$$
Since $(A,\lambda,v)\in \W$, we deduce that $\dot\eta=\dot\lambda$ and $\Pi_{v^\perp}\dot w=\dot v$. Then, 
$$
 \left(DF_A(\lambda,v)|_{\K\times {v'}^{\perp}}\right)^{-1} \, DF_A(\lambda,v)|_{\K\times v^{\perp}}(\dot\lambda,\dot v)=\big(\dot\lambda,({\Pi_{v^{\perp}}|}_{{v'}^{\perp}})^{-1}(\dot v)\big).
$$
(ii): Taking the canonical norm of $\K\times\K^n$ in (i), and maximizing on the unit sphere in $\K\times v^{\perp}\subset\K\times\K^n$, the assertion (ii) follows from Lemma~\ref{lem:tecvvpr}.
\\
(iii): Note that
\begin{align*}
\lefteqn{\|\left(DF_A(\lambda,v)|_{\K\times {v'}^{\perp}}\right)^{-1}\|  \leq }
\\
&&\|\left(DF_A(\lambda,v)|_{\K\times {v'}^{\perp}}\right)^{-1}\,
DF_A(\lambda,v)|_{\K\times v^{\perp}}\|\,\|\left(DF_A(\lambda,v)|_{\K\times
v^{\perp}}\right)^{-1}\|,
\end{align*}
then the result follows from (ii).
\end{proof}

\begin{lem}\label{lem:derseg}
Let $A\in\mnk$, and $(\lambda,v)\in\K\times\K^n$. Then $\|D^2F_A(\lambda,v)\|\leq 1$.
\end{lem}
\begin{proof}
 Differentiating $F_A$ twice, we get
$$D^2F_A(\lambda,v)(\dot\lambda,\dot v)(\dot\eta,\dot u)=\dot\lambda \dot u+ \dot\eta\dot v,\quad\mbox{for all}\quad (\dot\lambda,\dot v),\,(\dot\eta,\dot u)\in\K\times \K^n. $$
Then, 
\begin{align*}
\|D^2F_A(\lambda,v)(\dot\lambda,\dot v)(\dot\eta,\dot u)\| &\leq |\dot\lambda|\, \|\dot u\|+ \|\dot
v\|\,|\dot\eta|\\\
& \leq  (|\dot\lambda|^2+\|\dot v\|^2)^{1/2}\, (|\dot\eta|^2+\|\dot u\|^2)^{1/2},
\end{align*} 
where the second inequality follows from Cauchy-Schwarz inequality.
\end{proof}

Recall Neumann's series result (see for example Stewart--Sun~\cite{St-S}):
\begin{lem}\label{lem:NL}
 Let $\mathbb{E}$ be a Hermitian space, and $A, I_\mathbb{E}\,:\mathbb{E}\to \mathbb{E}$ be linear operators where $I_\mathbb{E}$ is the identity. If $\|A-I_\mathbb{E}\|<1$, then $A$ is invertible and
$$
\|A^{-1}\|\leq\frac{1}{1-\|A-I_\mathbb{E}\|}.
\qed
$$
\end{lem}

\begin{prop}\label{prop:NewtonAfindt}
Let $0<c\leq 1/(2\sqrt{2})$. \\
Let $A\in\mnk$ and $(\lambda,v)\in\K\times\K^n$, such that $\|v\|=1$ and $(A,\lambda,v)\in\W$. Let $(\lambda_0,v_0)\in\K\times\prk$.
If
$$(|\lambda_0-\lambda|^2+d_T(v_0,v)^2)^{1/2}<\frac{c}{
\left\|(DF_{A}\left(\lambda,v\right)|_{\K\times
v^{\perp}})^{-1} \right\| },$$
 then the sequence $(\lambda_k,v_k):=N^k_A(\lambda_0,v_0)$ satisfies
$$
(|\lambda_k-\lambda|^2+d_T(v_k,v)^2)^{1/2}\leq
  \sqrt{2}\,\delta(\sqrt{2}\,c)\,\left(\frac{1}{2}\right)^{2^k-1}\,
(|\lambda_0-\lambda|^2+d_T(v_0,v)^2)^{1/2},
$$
for all positive integers $k$, where  $\delta(c):=c/(1-c)$.
\end{prop}
\begin{proof}
Take a representative of ${v_0}$ such that $\pes{v-{v_0}}{v_0}=0$. Thus we have
${\|v_0\|\, d_T(v,{v_0})=\|v-{v_0}\|}$ and $\|v_0\|\leq 1$.

In particular, the hypothesis implies that
$$
\|(DF_A(\lambda,v)|_{\K\times
v^\perp})^{-1}\|\,\|({\lambda_0}-\lambda,{v_0}-v)\|< c.
$$

Taylor's expansion of $F_A$ and $DF_A$ in a neighbourhood of $(\lambda,v)$ are given by
\begin{equation}\label{eq:taylor}
F_A({\lambda'},{v'}) =DF_A(\lambda,v)({\lambda'}-\lambda,{v'}-v) + \frac12\, D^2F_A(\lambda,v)({\lambda'}-\lambda,{v'}-v)^2,
\end{equation}
and
$$
DF_A({\lambda'},{v'})=DF_A(\lambda,v) + D^2F_A(\lambda,v)({\lambda'}-\lambda,{v'}-v).
$$
One has
\begin{align*}
\lefteqn{\left(DF_A(\lambda,v)\big|_{\K\times {v_0}^\perp}\right)^{-1}\,
DF_A({\lambda_0},{v_0})\big|_{\K\times {v_0}^\perp} -
I_{\K\times {v_0}^\perp} = }\\
& \qquad=
\left(DF_A(\lambda,v)\big|_{\K\times
{v_0}^\perp}\right)^{-1}
\,\left( DF_A({\lambda_0},{v_0})\big|_{\K\times {v_0}^\perp}- DF_A(\lambda,v)\big|_{\K\times {v_0}^\perp} \right) \\
& \qquad= \left(DF_A(\lambda,v)\big|_{\K\times
{v_0}^\perp}\right)^{-1} \,
D^2F_A(\lambda,v))({\lambda_0}-\lambda,{v_0}-v)\big|_{\K\times
{v_0}^\perp} .
 \end{align*}
Then, taking norms, we
get
\begin{align*}
&\left\|\left(DF_A(\lambda,v)\big|_{\K\times
{v_0}^\perp}\right)^{-1}\, DF_A({\lambda_0},{v_0})\big|_{\K\times
{v_0}^\perp} -
I_{\K\times {v_0}^\perp}\right\| \leq \\
& \qquad \qquad \qquad\leq    \left\| (DF_A(\lambda,v)|_{\K\times
{v_0}^\perp})^{-1}\right\|
\, \left\| D^2F_A(\lambda,v))({\lambda_0}-\lambda,{v_0}-v)\right\|\\
& \qquad \qquad \qquad\leq 
\frac{1}{\cos(d_{\mathbb{P}}(v,{v_0}))}\, \left\|(DF_A(\lambda,v)|_{\K\times
{v}^\perp})^{-1}\right\| \,
\|({\lambda_0}-\lambda,{v_0}-v)\|, 
\end{align*}
where the last
inequality follows from {Lemma~\ref{lem:tec}} and
{Lemma~\ref{lem:derseg}}.
\\
In the range of angles under consideration we have $\|v_0\|=\cos(d_{\mathbb{P}}(v,{v_0}))\geq 1/\sqrt{2}$.
Then, by the condition $0<c\leq 1/(2\sqrt{2})$, we can deduce from Lemma
\ref{lem:NL}  that $ DF_A({\lambda_0},{v_0})\big|_{\K\times {v_0}^\perp}$ is
invertible and
\begin{align}
\label{eq:normin{v_0}}
&\left\|\left(DF_A({\lambda_0},{v_0})\big|_{\K\times {v_0}^\perp}\right)^{-1}\, DF_A(\lambda,v)\big|_{\K\times {v_0}^\perp} \right\|\leq\\
&\qquad\leq \frac{1}{1-\frac{1}{\cos(d_{\mathbb{P}}(v,{v_0}))}\, \left\|
(DF_A(\lambda,v)|_{\K\times v^\perp})^{-1}\right\| \,
\|({\lambda_0}-\lambda,{v_0}-v)\|}.\nonumber
\end{align}

Furthermore,
\begin{align*} 
&N_A({\lambda_0},{v_0})-(\lambda,v)= \\
&\qquad\qquad=  ({\lambda_0}-\lambda,{v_0}-v)- \left(DF_A({\lambda_0},{v_0})\big|_{\K\times {v_0}^\perp}\right)^{-1}\, F_A({\lambda_0},{v_0})\\
&\qquad\qquad  = \left(DF_A({\lambda_0},{v_0})\big|_{\K\times {v_0}^\perp}\right)^{-1}\, \\
 &\quad\qquad \qquad \qquad  \,\left( DF_A({\lambda_0},{v_0})\big|_{\K\times {v_0}^\perp}({\lambda_0}-\lambda,{v_0}-v)-F_A({\lambda_0},{v_0})\right).
\end{align*}
Then, from (\ref{eq:taylor})  we get
\begin{align*} 
& N_A({\lambda_0},{v_0})-(\lambda,v) = \\
&\qquad\qquad = \frac{1}{2}\,\left(DF_A({\lambda_0},{v_0})\big|_{\K\times
{v_0}^\perp}\right)^{-1}\,
D^2F_A(\lambda,v)({\lambda_0}-\lambda,{v_0}-v)^2.
\end{align*}
 Taking the canonical norm in $\K\times\K^n$, we get
\begin{align*} 
&\|N_A({\lambda_0},{v_0})-(\lambda,v)\| \leq \\
 &\quad \leq
\frac{1}{2}\, \left\| (DF_A({\lambda_0},{v_0})|_{\K\times
{v_0}^\perp})^{-1}\right\| \, \|
D^2F_A(\lambda,v)({\lambda_0}-\lambda,{v_0}-v)^2\|.
\end{align*}
Then,  from (\ref{eq:normin{v_0}})  and {Lemma~\ref{lem:tec}},
\begin{multline}\label{eq:distvunovcero}
{\|N_A({\lambda_0},{v_0})-(\lambda,v)\|\leq }\\
\leq \frac{\sqrt{2}\, \left\| (DF_A(\lambda,v)|_{\K\times
{v}^\perp})^{-1}\right\| \,\frac{1}{2}\,
 \|D^2F_A(\lambda,v)({\lambda_0}-\lambda,{v_0}-v)^2\|.}{1-\sqrt{2}\,
\left\|(DF_A(\lambda,v)|_{\K\times
{v}^\perp})^{-1}\right\| \,
\|({\lambda_0}-\lambda,{v_0}-v)\|}
\end{multline}
Therefore Lemma~\ref{lem:derseg} yields
\begin{multline*}
{\|N_A({\lambda_0},{v_0})-(\lambda,v)\|\leq }\\
\leq \frac{\sqrt{2}\, \left\| (DF_A(\lambda,v)|_{\K\times
{v}^\perp})^{-1}\right\| \,
\|({\lambda_0}-\lambda,{v_0}-v)\|}{1-\sqrt{2}\, 
\left\|(DF_A(\lambda,v)|_{\K\times {v}^\perp})^{-1}\right\| \,
\|({\lambda_0}-\lambda,{v_0}-v)\|}\,
\frac{1}{2}\,\|({\lambda_0}-\lambda,{v_0}-v)\|.\nonumber
\end{multline*}
Then,
\begin{multline*}
\lefteqn{\|N_A({\lambda_0},{v_0})-(\lambda,v)\|\leq }\\
\leq \frac{\sqrt{2}\, \left\| (DF_A(\lambda,v)|_{\K\times
{v}^\perp})^{-1}\right\| \,
(|{\lambda_0}-\lambda|^2+d_T({v_0},v)^2)^{1/2}}{1-\sqrt{2}\, 
\left\|(DF_A(\lambda,v)|_{\K\times {v}^\perp})^{-1}\right\| \,
(|{\lambda_0}-\lambda|^2+d_T({v_0},v)^2)^{1/2}}\, \\
\,
\frac{1}{2}\,(|{\lambda_0}-\lambda|^2+d_T({v_0},v)^2)^{1/2}.\nonumber
\end{multline*}

Let $(\lambda_1,v_1):=N_A(\lambda_0,v_0)$.\\
From the proof of {Lemma~\ref{lem:derseg}} we have
$D^2F_A(\lambda,v)(\lambda_0-\lambda,v_0-v)^2=2(\lambda_0-\lambda)(v_0-v)$,
then, from (\ref{eq:distvunovcero}) one can deduce that
$\|v_1-v\|<\delta(\sqrt{2}\,c)\|v_0-v\|$, where ${\delta(c)=c/(1-c)}$. Since
$c\leq 1/(2\sqrt{2})$, we have
$\delta(\sqrt{2}\,c)\leq 1$, then from Lemma 2--(4) of Blum et al.~\cite[pp. 264]{B-C-S-S} we get
$$
d_T(v_1,v)\leq
\frac{\|v_1-v\|}{\|v_0\|}\leq \sqrt{2}\,\|v_1-v\|.
$$
Hence
\begin{multline}\label{eq:induno}
\lefteqn{( |\lambda_1-\lambda|^2+ d_T(v_1,v)^2  )^{1/2}\leq} \\
\leq \frac{2\, \left\| (DF_A(\lambda,v)|_{\K\times
{v}^\perp})^{-1}\right\| \,
(|\lambda_0-\lambda|^2+d_T(v_0,v)^2)^{1/2}}{1-\sqrt{2}\, 
\left\|(DF_A(\lambda,v)|_{\K\times {v}^\perp})^{-1}\right\| \,
(|\lambda_0-\lambda|^2+d_T(v_0,v)^2)^{1/2}}
\,\\
\,\frac{1}{2}\,(|\lambda_0-\lambda|^2+d_T(v_0,v)^2)^{1/2}.
\end{multline}

Therefore
\beq\label{eq:preinduccion}
( |\lambda_1-\lambda|^2+ d_T(v_1,v)^2  )^{1/2}\leq
\sqrt{2}\,\delta(\sqrt{2}\,c)\,
\frac12\,(|\lambda_0-\lambda|^2+d_T(v_0,v)^2)^{1/2}.
\eeq

From (\ref{eq:preinduccion}), (\ref{eq:induno}), and the fact that
$\delta(\sqrt{2}\,c)\leq 1$, working by induction we get
$$
(|\lambda_k-\lambda|^2+d_T(v_k,v)^2)^{1/2}\leq
  \sqrt{2}\,\delta(\sqrt{2}\,c)\,\left(\frac{1}{2}\right)^{2^k-1}\,
(|\lambda_0-\lambda|^2+d_T(v_0,v)^2)^{1/2},
$$
for all $k>0$, where $(\lambda_k,v_k):=N^k_A(\lambda_0,v_0)$ .
\end{proof}

\begin{prop}\label{prop:muvsgamma} Let $(A,\lambda,v)\in\W$, such that $\|A\|_F=1$ and $\|v\|=1$. Then,
 $$
\mu(A,\lambda,v)\leq \|(DF_{A}\left(\lambda,v\right)|_{\K\times
v^{\perp}})^{-1} \|  \leq 2 \, \mu(A,\lambda,v).
$$
\end{prop}
\begin{proof}
 Since the action of $\mathbb{U}_n(\K)$ on $\prk$ is transitive, by Remark~\ref{rem:unitinv}, we may assume that $v=e_1$, where $e_1,\ldots,e_n$ is the canonical basis of $\K^n$. Then in this basis we have.
$$
A=
\begin{pmatrix}
 \lambda & a\\
0 & \hat A
\end{pmatrix}, 
$$
where $a\in\K^{1\times(n-1)}$, and $\hat{A}=\Pi_{{e_1}^\perp}A|_{{e_1}^\perp}\in\K^{(n-1)\times(n-1)}$. 

Recall that $DF_A(\lambda,{e_1})(\dot\lambda,\dot {v})=\dot\lambda {e_1} +(\lambda I_n-A)\dot {v}$. Then in the basis $(1,0)$, $(0,e_2),\ldots, (0,e_n)$ of $\K\times {e_1}^\perp$ we have
$$
DF_A(\lambda,{e_1})|_{\K\times {e_1}^\perp}=
\begin{pmatrix}
 1 & -a\\
0 &  \lambda I_{n-1}-\hat{A}
\end{pmatrix}.
$$

Note that $(DF_A(\lambda,{e_1})|_{\K\times
{e_1}^\perp})^{-1}=
\begin{pmatrix}
 1 & a  (\lambda I_{n-1}-\hat{A})^{-1}\\
0 &  (\lambda I_{n-1}-\hat{A})^{-1}
\end{pmatrix}.
$
Hence $$\|(DF_{A}\left(\lambda,{e_1}\right)|_{\K\times
{e_1}^{\perp}})^{-1} \|
\geq \max\{1, \| (\lambda I_{n-1}-\hat{A})^{-1} \|\}=\mu(A,\lambda,{e_1}).$$
On the other hand,
\begin{align*}
& \|(DF_{A}\left(\lambda,{e_1}\right)|_{\K\times {e_1}^{\perp}})^{-1} \|\leq \\
&\qquad\leq \left\|
\begin{pmatrix}
 1 & a  (\lambda I_{n-1}-\hat{A})^{-1}\\
0 &  0
\end{pmatrix}
\right\|+ \left\|
\begin{pmatrix}
 0 & 0\\
0 &  (\lambda I_{n-1}-\hat{A})^{-1}
\end{pmatrix}
\right\|\\
&\qquad \leq
\max\{1,\|a  (\lambda I_{n-1}-\hat{A})^{-1}\|\}+  \|(\lambda I_{n-1}-\hat{A})^{-1}\| \\
&\qquad\leq 2\,\mu(A,\lambda,{e_1}),
\end{align*}
where the last inequality follows from ${\|a\|\leq\|A\|_F=1}$.
\end{proof}

\begin{rem}
In the last proposition the result may be not longer true when ${\|A\|_F=1}$ or $\|v\|=1$ are not satisfied.
\end{rem}

\subsubsection{Proof of Proposition~\ref{prop:NewtonAfin}}
\begin{proof}[Proof of {Proposition~\ref{prop:NewtonAfin}}]
Pick a representative of $v$ and $v_0$ such that $\|v\|=1$ and $\pes{v-v_0}{v_0}=0$. Then the proof follows directly from {Proposition~\ref{prop:NewtonAfindt}}, {Proposition~\ref{prop:muvsgamma}} and {Lemma~\ref{lem:reldist}}.
\end{proof}

\medskip

\subsection{Proof of Lemma~\ref{des:dt}}

\begin{lem}\label{lem:apendice}
Let $A\in\mnk$ such that $\|A\|_F=1$. Let $\lambda,\,\lambda'
\in\K$ such that $|\lambda|\leq 1$.
\begin{enumerate}
 \item If $|\lambda'-\lambda|\leq c$ for some $0\leq c < \sqrt{2}$, then there
exists $\beta_c\,>1$ such that
$$
d_\mathbb{P}((A,\lambda),(A,\lambda))\leq \beta_c\,  |\lambda'-\lambda|.
$$
One may choose $\beta_c=(1-c^2/2)^{-1/2}$.
\item
If $d_{\mathbb{P}}((A,\lambda),(A,\lambda'))\leq \hat\theta$ for some $0\leq
\hat\theta < \pi/4$, then there exist $R_{\hat\theta}>1$ such that
$$
|\lambda'-\lambda|\leq R_{\hat\theta} \,
d_{\mathbb{P}}((A,\lambda),(A,\lambda')).
$$
One may choose $R_{\hat\theta}=[\sqrt{2}/\cos(\hat\theta+\pi/4)^3]^{1/2}$.
\end{enumerate}
\end{lem}
\begin{proof}
Let $\theta:=d_\mathbb{P}((A,\lambda),(A,\lambda'))$. (Thus $0\leq\theta\leq\hat\theta$.)
By the law of cosines we know that
$$
|\lambda-\lambda'|^2 = 1 +|\lambda|^2 +1 +|\lambda'|^2
-2\,\sqrt{1 +|\lambda|^2}\,\sqrt{1 +|\lambda'|^2}\,\cos
\theta.
$$
Then,
\begin{eqnarray}
|\lambda-\lambda'|^2 &=& \left(\sqrt{1 +|\lambda|^2} -\sqrt{1
+|\lambda'|^2} \right)^2  + \label{eq:cosinelaw}\\
&& \qquad+ 2\, \sqrt{1 +|\lambda|^2}\, \sqrt{1
+|\lambda|^2}\, (1-\cos\theta).\nonumber
\end{eqnarray}
From (\ref{eq:cosinelaw}) we get that
$$
|\lambda-\lambda'|^2 \geq 2\, \sqrt{1 +|\lambda|^2}\,
\sqrt{1 +|\lambda|^2}\, (1-\cos\theta),
$$
i.e.,
\beq\label{eq:cosinelaw2}
1-\cos\theta \leq \frac{|\lambda-\lambda'|^2}{2\, \sqrt{1
+|\lambda|^2}\, \sqrt{1 +|\lambda|^2}}
\leq\frac{{|\lambda'-\lambda|}^2}{2}.
\eeq
Therefore $1-\cos\theta \leq \frac{c^2}{2}$, and hence the angle $\theta$ is
bounded above by $\arccos{(1-c^2/2)}$.
By the Taylor expansion of cosine near $0$ we get the bound
$$
\theta^2\leq \frac{2}{1-c^2/2}\, (1-\cos\theta).
$$
Then, from (\ref{eq:cosinelaw2}) we can deduce the upper bound in (1).

For the lower bound in (2), we rewrite the cosine law and get:
\begin{eqnarray*}
|\lambda-\lambda'|^2 & = &
\left(
\frac{|\lambda|^2-|\lambda'|^2}{\sqrt{1 +|\lambda|^2}+\sqrt{1
+|\lambda'|^2}}
\right)^2+\\ 
&&\qquad\qquad\qquad\qquad\qquad+
2\sqrt{1 +|\lambda|^2}\,\sqrt{1
+|\lambda'|^2}\,(1-\cos\theta).
\end{eqnarray*}
Since $\left|| \lambda|-|\lambda'|\right|\leq|\lambda-\lambda'|$ and
$1-\cos\theta\leq \theta^2/2$, then,
\begin{eqnarray}\label{eq:cosinelaw3}
&&|\lambda-\lambda'|^2 \leq
 \left(
\frac{|\lambda|+|\lambda'|}{\sqrt{1 +|\lambda|^2}+\sqrt{1
+|\lambda'|^2}}
\right)^2\,|\lambda-\lambda'|^2+
\\
&&\qquad\qquad\qquad\qquad\qquad+
\sqrt{1 +|\lambda|^2}\,\sqrt{1
+|\lambda'|^2}\,\theta^2\nonumber
\end{eqnarray}

Since $0\leq |\lambda|\leq 1$, it is easily seen that
$$
\frac{|\lambda|+|\lambda'|}{\sqrt{1 +|\lambda|^2}+\sqrt{1
+|\lambda'|^2}}
\leq
\frac{1+|\lambda'|}{\sqrt{2}+\sqrt{1 +|\lambda'|^2}}.
$$
Furthermore, by elementary arguments one can see that $|\lambda'|\leq
\tan(\hat\theta+\pi/4)$, and therefore
\begin{align*}
\frac{|\lambda|+|\lambda'|}{\sqrt{1 +|\lambda|^2}+\sqrt{1
+|\lambda'|^2}}
&\leq \frac{1+ \tan(\hat\theta+\pi/4)}{\sqrt2+\sqrt{1+\tan(\hat\theta+\pi/4)^2}} \\
&\leq  
\frac{ \tan(\hat\theta+\pi/4)}{\sqrt{1+\tan(\hat\theta+\pi/4)^2}}=
\sin(\hat\theta+\pi/4),
\end{align*}
where the second inequality holds since $\tan(\hat\theta+\pi/4)\geq1$.
Then, from (\ref{eq:cosinelaw3}),
\beqna
|\lambda-\lambda'|^2 \leq \frac{\sqrt{1
+|\lambda|^2}\,\sqrt{1
+|\lambda'|^2}}{\cos(\hat\theta+\pi/4)^2}\,\theta^2,
\eeqna
and hence
$$
|\lambda-\lambda'|^2 \leq
\frac{\sqrt{2}}{\cos(\hat\theta+\pi/4)^3}\,\theta^2.
$$
\end{proof}
\begin{rem}
 Note that if $(A,\lambda)\in{\pi_1}(\V)\subset \mathbb{P}(\mnk\times\K)$
then $|\lambda|\leq \|A\|_F$ is always satisfied.
\end{rem}

\subsubsection{Proof of Lemma~\ref{des:dt}}
\begin{proof}[Proof of Lemma~\ref{des:dt}]
 The proof of (1) and (2) follows directly from the definition of
$d_{\mathbb{P}^2}$ and {Lemma~\ref{lem:apendice}}.
\end{proof}

\end{document}